\documentclass[12pt]{article}

\usepackage{epic}
\usepackage{eepic}
\usepackage{epsf} 
\usepackage{ amssymb, amscd}
\usepackage{amsmath}
\usepackage{color}
\numberwithin{equation}{section}
\usepackage{epsfig}
\usepackage{graphicx}

\def\cC            {{\mathcal{C}}}

\def\cM            {{\mathcal{M}}}

\def\cO            {{\mathcal{O}}}
\def\cP            {{\mathcal{P}}}

\def\cV            {{\mathcal{V}}}

\def\bbC           {\mathbb{C}}

\def\bbH          {\mathbb{H}}
\def\bbM           {\mathbb{M}}

\def\bbQ           {\mathbb{Q}}
\def\bbR           {\mathbb{R}}

\def\bbZ           {\mathbb{Z}}

\def\si{\sigma}
\def\tsi{\widetilde{\sigma}}

\def\I{{\rm i}}
\def\sdprod{{\times\!\vrule height5pt depth0pt width0.4pt\,}}
\def\eps{\epsilon}

\title{Much ado about Mathieu}
\author{
 {\sc Terry  Gannon }\\ 
 {\footnotesize Department of Mathematics, University of Alberta,}\\
{\footnotesize Edmonton, Alberta, Canada T6G 2G1}\\
{\footnotesize e-mail: {\tt tgannon@math.ualberta.ca}} }

\begin{document}
\maketitle

\begin{abstract} Eguchi, 
Ooguri and Tachikawa have observed that the elliptic genus of type II string theory on K3 
surfaces appears to possess a Moonshine for the largest Mathieu group. Subsequent work by
several people established a candidate for the elliptic genus twisted by each element of $M_{24}$. 
In this paper we prove that the resulting sequence of class functions are true characters of
$M_{24}$, proving the Eguchi-Ooguri-Tachikawa conjecture. The integrality of
multiplicities is proved using a small generalisation of Sturm's Theorem, while positivity
involves a modification of a method of Hooley.  We also prove the evenness
property of the multiplicities, as conjectured by several authors. We also identify the role group cohomology
plays in both K3-Mathieu Moonshine and Monstrous Moonshine; in particular this gives
 a cohomological interpretation for the non-Fricke
elements in Norton's Generalised Monstrous Moonshine conjecture.  We investigate
the proposal of Gaberdiel-Hohenegger-Volpato that K3-Mathieu Moonshine lifts to
the Conway group Co$_1$. 
\end{abstract}

{\footnotesize
\tableofcontents
}

\section{Introduction}

The elliptic genus (a.k.a. partition function) of a nonlinear sigma model with K3 target space is a very special function. On general 
grounds, it is a weak Jacobi form of index one and weight zero for $\bbZ^2\sdprod\mathrm{SL}_2(\bbZ)$, and is therefore equal to the 
unique (up to scaling) such function, $\phi_{0,1}(z,\tau)$. This implies it is ÔtopologicalÕ in the sense that it is independent of where you are on the 20-(complex) dimensional  K3 moduli space,
or indeed of where you are on the 40-dimensional moduli space of K3 sigma models.

The world-sheet description of these string theories is as an $N = 4$ superconformal field
theory, and thus the whole state space can be organised as an $N = 4$ superconformal 
representation. In particular,  the elliptic genus can be decomposed into sums of  elliptic 
genera of $N = 4$ superconformal representations. Eguchi, 
Ooguri and Tachikawa \cite{EOT} observed that the multiplicities with which these $N = 4$ elliptic genera contribute 
to the elliptic genus are dimensions of $M_{24}$-group representations.

According to Conway, the sporadic simple group $M_{24}$ is the most remarkable of all finite groups (p.300, \cite{CN}), with a great 
wealth of structure and applications. For instance it is the automorphism group of the Golay code, a 
stabilizer in the Leech lattice, and all symplectic symmetries of K3 surfaces can be embedded
 in it (in fact in a 
subgroup $M_{23}$) \cite{Mu,Ko}. Its order (size) is $244\,823\,040$, and class number
(i.e. number of conjugacy classes) is 26.  The character  table is given in Table 1;  we  retain the order of irreps(= irreducible representations) in \cite{Atl}, relate with a bar the 
 complex conjugate irrep (when nonisomorphic), and write  $\rho_0=1$ for the trivial representation. We write there $\alpha_t^s=(s+t\I\sqrt{7})/2$, $\beta_t^s=(s+t\I\sqrt{15})/2$, and $\gamma_t^s=(s+t\I\sqrt{23})/2$, for all choices of  signs $s,t\in\{\pm 1\}$.

K3-Mathieu Moonshine was pushed further --- indeed, made well-defined ---
by the work of  Cheng \cite{Ch}, Gaberdiel, Hohenegger \& Volpato 
\cite{GHV1,GHV2} and Eguchi \& Hikami \cite{EH} who calculated the analogue $\phi_g(\tau,z)$
of the McKay-Thompson
series here. The resulting functions are weak Jacobi forms for $\bbZ^2\sdprod \Gamma_0(|g|)$ 
up to certain phases ($|g|$ denotes the order of $g\in M_{24}$). At the conclusion of this
paper (see also \cite{Ron})  we
 identify the cohomological source of these phases. The K3-Mathieu Moonshine fits beautifully
 into this general cohomological framework of conformal field theory (CFT) \cite{Ron}.

Once we know the weak Jacobi forms $\phi_g$, we obtain modular forms $f_g$ from
\eqref{phigfg} below and can read off the class functions $H_n$ from \eqref{eq:fg} and \eqref{dmz}
 (a class function is a function constant on conjugacy classes).
These $H_n$ can also be obtained explicitly from \eqref{Hk} below. In any case, 
our most important result is the proof that these class functions $H_n$  are all true representations of $M_{24}$. This can be
regarded as a proof of the weak form of the K3-Mathieu moonshine conjecture.

\medskip\noindent\textbf{Theorem A.} \textit{Each $H_n$ ($n\ge 1$) is a true character of $M_{24}$
(i.e. a sum of irreducible characters of $M_{24}$).}

\vfill\eject
\centerline {{\bf Table 1.} Character table  of $M_{24}$}
$${\scriptsize\vbox{\tabskip=0pt\offinterlineskip
  \def\tablerule{\noalign{\hrule}}
  \halign to 6.9in{
  \strut#& \vrule#&     
\hfil#&\hfil#&\hfil#&\hfil#&\hfil#&\hfil#&\hfil#&\hfil#&\hfil#&\hfil#&\hfil#&\hfil#&\hfil#&\hfil#&\hfil#&\hfil#&\hfil#&\hfil#&\hfil#&\hfil#&
\hfil#&\hfil#&\hfil#&\hfil#&\hfil#\tabskip=0pt\cr$K_g$&\,\,1A\,\,&\,\,2A\,\,&\,\,2B\,\,&\,\,3A\,\,&\,\,3B\,\,&\,\,4A\,\,&\,\,4B\,\,&\,\,4C\,\,&\,\,5A\,\,&\,\,6A\,\,&\,\,6B\,\,&\,\,7A\,\,&\,\,7B\,\,&\,\,8A\,\,&\,10A\,&\,11A\,&\,12A\,&\,12B\,&\,\,14A\,\,&\,\,14B\,\,&\,\,15A\,\,&\,\,15B\,\,&\,\,21A\,\,&\,\,21B\,\,&\,\,23A\,\,&\,\,23B\cr\tablerule  
$1$&\,\,1\,\,&\,\,1\,\,&\,\,1\,\,&\,\,1\,\,&\,\,1\,\,&\,\,1\,\,&\,\,1\,\,&\,\,1\,\,&\,\,1\,\,&\,\,1\,\,&\,\,1\,\,&\,\,1\,\,&\,\,1\,\,&\,\,1\,\,&\,\,1\,\,&\,\,1\,\,&\,\,1\,\,&\,\,1\,\,&\,\,1\,\,&\,\,1\,\,&\,\,1\,\,&\,\,1\,\,&\,\,1\,\,&\,\,1\,\,&\,\,1\,\,&\,\,1\cr $\rho_1$
&\,\,23\,\,&\,\,7\,\,&\,\,-1\,\,&\,\,5\,\,&\,\,-1\,\,&\,\,-1\,\,&\,\,3\,\,&\,\,-1\,\,&\,\,3\,\,&\,\,1\,\,&\,\,-1\,\,&\,\,2\,\,&\,\,2\,\,&\,\,1\,\,&\,\,-1\,\,&\,\,1\,\,&\,\,-1\,\,&\,\,-1\,\,&\,\,0\,\,&\,\,0\,\,&\,\,0\,\,&\,\,0\,\,&\,\,-1\,\,&\,\,-1\,\,&\,\,0\,\,&\,\,0\cr$\rho_2$
&\,\,45\,\,&\,\,-3\,\,&\,\,5\,\,&\,\,0\,\,&\,\,3\,\,&\,\,-3\,\,&\,\,1\,\,&\,\,1\,\,&\,\,0\,\,&\,\,0\,\,&\,\,-1\,\,&\,\,$\alpha_+^-$\,\,&\,\,$\alpha_-^-$\,\,&\,\,-1\,\,&\,\,0\,\,&\,\,1\,\,&\,\,0\,\,&\,\,1\,\,&\,\,$\alpha_-^+$\,\,&\,\,$\alpha_+^+$\,\,&\,\,0\,\,&\,\,0\,\,&\,\,$\alpha_+^-$\,\,&\,\,$\alpha_-^-$\,\,&\,\,-1\,\,&\,\,-1\cr$\overline{\rho_2}$
&\,\,45\,\,&\,\,-3\,\,&\,\,5\,\,&\,\,0\,\,&\,\,3\,\,&\,\,-3\,\,&\,\,1\,\,&\,\,1\,\,&\,\,0\,\,&\,\,0\,\,&\,\,-1\,\,&\,\,$\alpha_-^-$\,\,&\,\,$\alpha_+^-$\,\,&\,\,-1\,\,&\,\,0\,\,&\,\,1\,\,&\,\,0\,\,&\,\,1\,\,&\,\,$\alpha_+^+$\,\,&\,\,$\alpha_-^+$\,\,&\,\,0\,\,&\,\,0\,\,&\,\,$\alpha_-^-$\,\,&\,\,$\alpha_+^-$\,\,&\,\,-1\,\,&\,\,-1\cr$\rho_3$
&\,\,231\,\,&\,\,7\,\,&\,\,-9\,\,&\,\,-3\,\,&\,\,0\,\,&\,\,-1\,\,&\,\,-1\,\,&\,\,3\,\,&\,\,1\,\,&\,\,1\,\,&\,\,0\,\,&\,\,0\,\,&\,\,0\,\,&\,\,-1\,\,&\,\,1\,\,&\,\,0\,\,&\,\,-1\,\,&\,\,0\,\,&\,\,0\,\,&\,\,0\,\,&\,\,$\beta_{+}^-$\,\,&\,\,$\beta_{-}^-$\,\,&\,\,0\,\,&\,\,0\,\,&\,\,1\,\,&\,\,1\cr$\overline{\rho_3}$
&\,\,231\,\,&\,\,7\,\,&\,\,-9\,\,&\,\,-3\,\,&\,\,0\,\,&\,\,-1\,\,&\,\,-1\,\,&\,\,3\,\,&\,\,1\,\,&\,\,1\,\,&\,\,0\,\,&\,\,0\,\,&\,\,0\,\,&\,\,-1\,\,&\,\,1\,\,&\,\,0\,\,&\,\,-1\,\,&\,\,0\,\,&\,\,0\,\,&\,\,0\,\,&\,\,$\beta_{-}^-$\,\,&\,\,$\beta_{+}^-$\,\,&\,\,0\,\,&\,\,0\,\,&\,\,1\,\,&\,\,1\cr$\rho_4$
&\,\,252\,\,&\,\,28\,\,&\,\,12\,\,&\,\,9\,\,&\,\,0\,\,&\,\,4\,\,&\,\,4\,\,&\,\,0\,\,&\,\,2\,\,&\,\,1\,\,&\,\,0\,\,&\,\,0\,\,&\,\,0\,\,&\,\,0\,\,&\,\,2\,\,&\,\,-1\,\,&\,\,1\,\,&\,\,0\,\,&\,\,0\,\,&\,\,0\,\,&\,\,-1\,\,&\,\,-1\,\,&\,\,0\,\,&\,\,0\,\,&\,\,-1\,\,&\,\,-1\cr$\rho_5$
&\,\,253\,\,&\,\,13\,\,&\,\,-11\,\,&\,\,10\,\,&\,\,1\,\,&\,\,-3\,\,&\,\,1\,\,&\,\,1\,\,&\,\,3\,\,&\,\,-2\,\,&\,\,1\,\,&\,\,1\,\,&\,\,1\,\,&\,\,-1\,\,&\,\,-1\,\,&\,\,0\,\,&\,\,0\,\,&\,\,1\,\,&\,\,-1\,\,&\,\,-1\,\,&\,\,0\,\,&\,\,0\,\,&\,\,1\,\,&\,\,1\,\,&\,\,0\,\,&\,\,0\cr$\rho_6$
&\,\,483\,\,&\,\,35\,\,&\,\,3\,\,&\,\,6\,\,&\,\,0\,\,&\,\,3\,\,&\,\,3\,\,&\,\,3\,\,&\,\,-2\,\,&\,\,2\,\,&\,\,0\,\,&\,\,0\,\,&\,\,0\,\,&\,\,-1\,\,&\,\,-2\,\,&\,\,-1\,\,&\,\,0\,\,&\,\,0\,\,&\,\,0\,\,&\,\,0\,\,&\,\,1\,\,&\,\,1\,\,&\,\,0\,\,&\,\,0\,\,&\,\,0\,\,&\,\,0\cr$\rho_7$
&\,\,770\,\,&\,\,-14\,\,&\,\,10\,\,&\,\,5\,\,&\,\,-7\,\,&\,\,2\,\,&\,\,-2\,\,&\,\,-2\,\,&\,\,0\,\,&\,\,1\,\,&\,\,1\,\,&\,\,0\,\,&\,\,0\,\,&\,\,0\,\,&\,\,0\,\,&\,\,0\,\,&\,\,-1\,\,&\,\,1\,\,&\,\,0\,\,&\,\,0\,\,&\,\,0\,\,&\,\,0\,\,&\,\,0\,\,&\,\,0\,\,&\,\,$\gamma_{+}^-$\,\,&\,\,$\gamma_{-}^-$\cr$\overline{\rho_7}$
&\,\,770\,\,&\,\,-14\,\,&\,\,10\,\,&\,\,5\,\,&\,\,-7\,\,&\,\,2\,\,&\,\,-2\,\,&\,\,-2\,\,&\,\,0\,\,&\,\,1\,\,&\,\,1\,\,&\,\,0\,\,&\,\,0\,\,&\,\,0\,\,&\,\,0\,\,&\,\,0\,\,&\,\,-1\,\,&\,\,1\,\,&\,\,0\,\,&\,\,0\,\,&\,\,0\,\,&\,\,0\,\,&\,\,0\,\,&\,\,0\,\,&\,\,$\gamma_{-}^-$\,\,&\,\,$\gamma_{+}^-$\cr$\rho_8$
&\,\,990\,\,&\,\,-18\,\,&\,\,-10\,\,&\,\,0\,\,&\,\,3\,\,&\,\,6\,\,&\,\,2\,\,&\,\,-2\,\,&\,\,0\,\,&\,\,0\,\,&\,\,-1\,\,&\,\,$\alpha_+^-$\,\,&\,\,$\alpha_-^-$\,\,&\,\,0\,\,&\,\,0\,\,&\,\,0\,\,&\,\,0\,\,&\,\,1\,\,&\,\,$\alpha_+^-$\,\,&\,\,$\alpha_-^-$\,\,&\,\,0\,\,&\,\,0\,\,&\,\,$\alpha_+^-$\,\,&\,\,$\alpha_-^-$\,\,&\,\,1\,\,&\,\,1\cr$\overline{\rho_8}$
&\,\,990\,\,&\,\,-18\,\,&\,\,-10\,\,&\,\,0\,\,&\,\,3\,\,&\,\,6\,\,&\,\,2\,\,&\,\,-2\,\,&\,\,0\,\,&\,\,0\,\,&\,\,-1\,\,&\,\,$\alpha_-^-$\,\,&\,\,$\alpha_+^-$\,\,&\,\,0\,\,&\,\,0\,\,&\,\,0\,\,&\,\,0\,\,&\,\,1\,\,&\,\,$\alpha_-^-$\,\,&\,\,$\alpha_+^-$\,\,&\,\,0\,\,&\,\,0\,\,&\,\,$\alpha_-^-$\,\,&\,\,$\alpha_+^-$\,\,&\,\,1\,\,&\,\,1\cr$\rho_9$
&\,\,1035\,\,&\,\,27\,\,&\,\,35\,\,&\,\,0\,\,&\,\,6\,\,&\,\,3\,\,&\,\,-1\,\,&\,\,3\,\,&\,\,0\,\,&\,\,0\,\,&\,\,2\,\,&\,\,-1\,\,&\,\,-1\,\,&\,\,1\,\,&\,\,0\,\,&\,\,1\,\,&\,\,0\,\,&\,\,0\,\,&\,\,-1\,\,&\,\,-1\,\,&\,\,0\,\,&\,\,0\,\,&\,\,-1\,\,&\,\,-1\,\,&\,\,0\,\,&\,\,0\cr$\rho_{10}$
&\,\,1035\,\,&\,\,-21\,\,&\,\,-5\,\,&\,\,0\,\,&\,\,-3\,\,&\,\,3\,\,&\,\,3\,\,&\,\,-1\,\,&\,\,0\,\,&\,\,0\,\,&\,\,1\,\,&\,\,$2\alpha_+^-$\,\,&\,\,$2\alpha_-^-$\,\,&\,\,-1\,\,&\,\,0\,\,&\,\,1\,\,&\,\,0\,\,&\,\,-1\,\,&\,\,0\,\,&\,\,0\,\,&\,\,0\,\,&\,\,0\,\,&\,\,$\alpha_-^+$\,\,&\,\,$\alpha_+^+$\,\,&\,\,0\,\,&\,\,0\cr$\overline{\rho_{10}}$
&\,\,1035\,\,&\,\,-21\,\,&\,\,-5\,\,&\,\,0\,\,&\,\,-3\,\,&\,\,3\,\,&\,\,3\,\,&\,\,-1\,\,&\,\,0\,\,&\,\,0\,\,&\,\,1\,\,&\,\,$2\alpha_-^-$\,\,&\,\,$2\alpha_+^-$\,\,&\,\,-1\,\,&\,\,0\,\,&\,\,1\,\,&\,\,0\,\,&\,\,-1\,\,&\,\,0\,\,&\,\,0\,\,&\,\,0\,\,&\,\,0\,\,&\,\,$\alpha_+^+$\,\,&\,\,$\alpha_-^+$\,\,&\,\,0\,\,&\,\,0\cr$\rho_{11}$
&\,\,1265\,\,&\,\,49\,\,&\,\,-15\,\,&\,\,5\,\,&\,\,8\,\,&\,\,-7\,\,&\,\,1\,\,&\,\,-3\,\,&\,\,0\,\,&\,\,1\,\,&\,\,0\,\,&\,\,-2\,\,&\,\,-2\,\,&\,\,1\,\,&\,\,0\,\,&\,\,0\,\,&\,\,-1\,\,&\,\,0\,\,&\,\,0\,\,&\,\,0\,\,&\,\,0\,\,&\,\,0\,\,&\,\,1\,\,&\,\,1\,\,&\,\,0\,\,&\,\,0\cr$\rho_{12}$
&\,\,1771\,\,&\,\,-21\,\,&\,\,11\,\,&\,\,16\,\,&\,\,7\,\,&\,\,3\,\,&\,\,-5\,\,&\,\,-1\,\,&\,\,1\,\,&\,\,0\,\,&\,\,-1\,\,&\,\,0\,\,&\,\,0\,\,&\,\,-1\,\,&\,\,1\,\,&\,\,0\,\,&\,\,0\,\,&\,\,-1\,\,&\,\,0\,\,&\,\,0\,\,&\,\,1\,\,&\,\,1\,\,&\,\,0\,\,&\,\,0\,\,&\,\,0\,\,&\,\,0\cr$\rho_{13}$
&\,\,2024\,\,&\,\,8\,\,&\,\,24\,\,&\,\,-1\,\,&\,\,8\,\,&\,\,8\,\,&\,\,0\,\,&\,\,0\,\,&\,\,-1\,\,&\,\,-1\,\,&\,\,0\,\,&\,\,1\,\,&\,\,1\,\,&\,\,0\,\,&\,\,-1\,\,&\,\,0\,\,&\,\,-1\,\,&\,\,0\,\,&\,\,1\,\,&\,\,1\,\,&\,\,-1\,\,&\,\,-1\,\,&\,\,1\,\,&\,\,1\,\,&\,\,0\,\,&\,\,0\cr$\rho_{14}$
&\,\,2277\,\,&\,\,21\,\,&\,\,-19\,\,&\,\,0\,\,&\,\,6\,\,&\,\,-3\,\,&\,\,1\,\,&\,\,-3\,\,&\,\,-3\,\,&\,\,0\,\,&\,\,2\,\,&\,\,2\,\,&\,\,2\,\,&\,\,-1\,\,&\,\,1\,\,&\,\,0\,\,&\,\,0\,\,&\,\,0\,\,&\,\,0\,\,&\,\,0\,\,&\,\,0\,\,&\,\,0\,\,&\,\,-1\,\,&\,\,-1\,\,&\,\,0\,\,&\,\,0\cr$\rho_{15}$
&\,\,3312\,\,&\,\,48\,\,&\,\,16\,\,&\,\,0\,\,&\,\,-6\,\,&\,\,0\,\,&\,\,0\,\,&\,\,0\,\,&\,\,-3\,\,&\,\,0\,\,&\,\,-2\,\,&\,\,1\,\,&\,\,1\,\,&\,\,0\,\,&\,\,1\,\,&\,\,1\,\,&\,\,0\,\,&\,\,0\,\,&\,\,-1\,\,&\,\,-1\,\,&\,\,0\,\,&\,\,0\,\,&\,\,1\,\,&\,\,1\,\,&\,\,0\,\,&\,\,0\cr$\rho_{16}$
&\,\,3520\,\,&\,\,64\,\,&\,\,0\,\,&\,\,10\,\,&\,\,-8\,\,&\,\,0\,\,&\,\,0\,\,&\,\,0\,\,&\,\,0\,\,&\,\,-2\,\,&\,\,0\,\,&\,\,-1\,\,&\,\,-1\,\,&\,\,0\,\,&\,\,0\,\,&\,\,0\,\,&\,\,0\,\,&\,\,0\,\,&\,\,1\,\,&\,\,1\,\,&\,\,0\,\,&\,\,0\,\,&\,\,-1\,\,&\,\,-1\,\,&\,\,1\,\,&\,\,1\cr$\rho_{17}$
&\,\,5313\,\,&\,\,49\,\,&\,\,9\,\,&\,\,-15\,\,&\,\,0\,\,&\,\,1\,\,&\,\,-3\,\,&\,\,-3\,\,&\,\,3\,\,&\,\,1\,\,&\,\,0\,\,&\,\,0\,\,&\,\,0\,\,&\,\,-1\,\,&\,\,-1\,\,&\,\,0\,\,&\,\,1\,\,&\,\,0\,\,&\,\,0\,\,&\,\,0\,\,&\,\,0\,\,&\,\,0\,\,&\,\,0\,\,&\,\,0\,\,&\,\,0\,\,&\,\,0\cr$\rho_{18}$
&\,\,5544\,\,&\,\,-56\,\,&\,\,24\,\,&\,\,9\,\,&\,\,0\,\,&\,\,-8\,\,&\,\,0\,\,&\,\,0\,\,&\,\,-1\,\,&\,\,1\,\,&\,\,0\,\,&\,\,0\,\,&\,\,0\,\,&\,\,0\,\,&\,\,-1\,\,&\,\,0\,\,&\,\,1\,\,&\,\,0\,\,&\,\,0\,\,&\,\,0\,\,&\,\,-1\,\,&\,\,-1\,\,&\,\,0\,\,&\,\,0\,\,&\,\,1\,\,&\,\,1\cr$\rho_{19}$
&\,\,5796\,\,&\,\,-28\,\,&\,\,36\,\,&\,\,-9\,\,&\,\,0\,\,&\,\,-4\,\,&\,\,4\,\,&\,\,0\,\,&\,\,1\,\,&\,\,-1\,\,&\,\,0\,\,&\,\,0\,\,&\,\,0\,\,&\,\,0\,\,&\,\,1\,\,&\,\,-1\,\,&\,\,-1\,\,&\,\,0\,\,&\,\,0\,\,&\,\,0\,\,&\,\,1\,\,&\,\,1\,\,&\,\,0\,\,&\,\,0\,\,&\,\,0\,\,&\,\,0\cr$\rho_{20}$
&\,\,10395\,\,&\,\,-21\,\,&\,\,-45\,\,&\,\,0\,\,&\,\,0\,\,&\,\,3\,\,&\,\,-1\,\,&\,\,3\,\,&\,\,0\,\,&\,\,0\,\,&\,\,0\,\,&\,\,0\,\,&\,\,0\,\,&\,\,1\,\,&\,\,0\,\,&\,\,0\,\,&\,\,0\,\,&\,\,0\,\,&\,\,0\,\,&\,\,0\,\,&\,\,0\,\,&\,\,0\,\,&\,\,0\,\,&\,\,0\,\,&\,\,-1\,\,&\,\,-1\cr
\noalign{\smallskip}}}}$$

The \textit{strong} conjecture says that these $\phi_g$  are the twisted elliptic genera of
an $N=4$ vertex operator superalgebra --- we have nothing to say about this. One would have 
liked to prove Theorem A by explicitly constructing these $M_{24}$-representations $H_n$.
This remains an important challenge. Instead, we prove Theorem A in two steps.

The first step in showing this is to prove that the $H_n$ are \textit{virtual characters},
i.e. linear combinations over the integers of irreducible characters.
It can be easily shown that the coefficients $H_n(g)$ will be \textit{integer-valued} class
functions, but this only
implies the class function is a linear combination over $\bbQ$ of irreducibles, with
denominators dividing the order $|G|$ of the group (244 million, in our case!). We prove
the $H_n$ are virtual  by verifying that for each $n$ the quantities $H_n(g)$ satisfy certain
congruences; we can verify this simultaneously for all $n$ by studying mod $p^\alpha$ reductions
of associated modular forms. (Our result Lemma 3 on mod $n$ reductions of modular
forms is a refinement of older results in the literature.)

The other and more difficult step in proving Theorem A, is to show that each $H_n(g)$ is a
\textit{nonnegative} (real) linear combination of irreducible characters. It is easy to reduce
this to showing that for any $g\ne 1$, $|H_n(g)|$ is small compared to $H_n(1)$ for all
$n>N$ for some  $N$ (the coefficients for $n\le N$ are then checked explicitly). 
Thanks to \cite{CD} we know 
\begin{equation}
H_n(g)=\cO\left(\frac{e^{\pi\sqrt{8n-1}/2|g|}}{\sqrt{8n-1}}\right)\end{equation}
where $|g|$ denotes the order of $g$, so for all  $n>N$, where $N$ is sufficiently large, $H_n(1)$ will dwarf $H_n(g)$
for $g\ne 1$, but we need to make this \textit{effective} by finding that $N$, and this is the hard
part of this positivity argument. It is tempting to guess we could follow Rademacher's calculation of 
effective bounds for the partition numbers, which looks like it should be similar, but
Rademacher's series was absolutely convergent whereas those in \cite{CD} aren't, so
the argument is much more delicate. But the analogy with partition numbers is still
useful:  \textit{generalised} Kloosterman sums arise here, and the corresponding zeta
functions are what we need to bound; but
 Whiteman \cite{Wh} showed long ago that the classical Kloosterman sum  has an elegant
 expression as a sparse sum, and Hooley \cite{Hoo} explained how to use the theory of
 binary quadratic forms to bound  series associated to similar sums, so we follow their lead.
 This involves though finding significant generalisations of Whiteman's and Hooley's
 results. 
 As an aside, our inequalities yield an independent proof of the convergence of the Rademacher sum
 expressions in \cite{CD}.

Integrality of these multiplicities is crucial. In contrast,  the significance of
positivity is not as  clear to this author --- after all, elliptic genus is a \textit{signed} trace. 
What happens here
should be contrasted with the $N=4$ $c=6$ toroidal theories, where the elliptic genus vanishes.
Positivity
could be a consequence of the minimality of this string theory (see e.g. the discussions
in \cite{Oo},\cite{We}) --- we return to this important question in the conclusion. 
In any case, it is intriguing to note that
the elliptic genera of the $M_{24}$-twisted modules (these have been recently obtained
in \cite{Ron}) also appear to be true representations (of the appropriate central extension of
centralisers in $M_{24}$).

Our methods of proving weak Moonshine are quite robust. In particular, they can surely be
applied to the twisted twining Mathieu Moonshine elliptic genera of \cite{Ron}, the `Umbral Moonshine' of \cite{CDH} and the PSL$_2(\bbZ_{11})$ $N=2$ Moonshine 
of \cite{EHb} (as well of course to Monstrous Moonshine itself).

It is a consequence of our proof that the multiplicities mult$_{\rho_i}(H_n)$ tend to infinity
with $n$, for each $i$, and all mult$_{\rho_i}(H_n)$ are strictly positive for $n\ge 25$. See Theorem 4 below.

It is common in the Mathieu Moonshine literature to emphasise the presence of \textit{mock} modular
forms. Indeed, the sums $q^{-1/8}\sum_{n=0}^\infty H_n(g)q^n$ are mock modular for each $g$,
not (usually) modular. We had no use however of mock modularity in this paper; however the true
modularity of the derived functions we call $f_g$ plays a crucial role.

In Section 3 we also obtain an evenness result conjectured by several authors:

\medskip\noindent\textbf{Theorem B.} \textit{Each head character $H_n$ is a sum of}
\begin{eqnarray}&&\{2, 2\rho_1,\rho_2+\overline{\rho_2}, \rho_3+\overline{\rho_3},2\rho_4,2\rho_5,2\rho_6,\rho_7+\overline{\rho_7},\rho_8+
\overline{\rho_8},2\rho_9,\nonumber\\
&&\qquad\quad\,\,\rho_{10}+\overline{\rho_{10}},2\rho_{11},2\rho_{12},2\rho_{13},2\rho_{14},2\rho_{15},2\rho_{16},2\rho_{17},2\rho_{18},2\rho_{19},2\rho_{20}\}\nonumber\end{eqnarray}

Theorems A and B are both used in \cite{CHM} to prove a conjecture in Umbral Moonshine
\cite{CDH}. Indeed, their Theorem 1.2 is far stronger than Conjecture 5.11 in \cite{CDH}
(specialised to $M_{24}$),
and a little weaker than Conjecture 5.12 in \cite{CDH}.
 Incidentally, it is curious that imaginary quadratic fields 
play a role in Umbral Moonshine (see Section 5.4 of \cite{CDH}) while the equivalent
theory of positive-definite quadratic forms plays a crucial role in our Section 4 proof.
Perhaps this can supply a deeper explanation for Conjecture 5.12 than would arise
from arguments as in \cite{CHM}.

Theorems A and B also imply that the elliptic genus of Enriques surfaces can be decomposed
into true characters of $M_{12}$ \cite{EHc} (their elliptic genus is half that of K3 surfaces,
so positivity comes from Theorem A and integrality from Theorem B). The VOA or string
theoretic interpretation of this observation seems quite obscure. 

 \cite{GHV3} made the intriguing proposal that in fact the much larger
 Conway groups Co$_0$ or  Co$_1$ may act. We prove (Gerald H\"ohn notified us that he also has a proof):
 
\medskip\noindent\textbf{Theorem C.} \textit{All $H_n$, as well as $H_{00}$, are
restrictions of virtual representations of Co$_0$.}\medskip

This was a nontrivial obstruction to their suggestion. However, as we explain below, the virtual
Co$_0$-representations are much larger than the underlying $M_{24}$-representations
(i.e. in the restriction a large  representation with negative multiplicity cancels a large representation with positive multiplicity,
leaving the small remainder $H_n$). Indeed, it now seems few people expect Mathieu Moonshine
to extend to Co$_0$. Of course Theorem C implies that the same conclusion necessarily 
holds when Co$_0$ is replaced there with any subgroup between $M_{24}$ and Co$_0$.
In the conclusion, we mention the other group that should be considered 
 a candidate for enhanced symmetry in Mathieu Moonshine.

In our paper, and indeed in much work on the subject, no relation involving K3 is used or
obtained. It seems possible to us that the connection with K3 may be illusory. This would be
very disappointing. We return to this a little more in the concluding section.

Because a fairly wide range of readers are potentially interested in parts of this paper, and
may wish to extend it to other Moonshines, we
have endeavoured to keep it as accessible and as detailed as possible. 

\section{K3-Mathieu Moonshine: Review}

Throughout this paper, write $e(x)=\exp(2\pi\I x)$. For readability we will often collapse
$a\equiv b$ (mod $m$) to $a\equiv_m b$. We will write $|g|$ for the order $|\langle g\rangle|$
of the element $g$ of a group. An excellent review of superconformal field theories and their
moduli spaces and elliptic genera, including $N=4$ $c=6$, is included in \cite{We}.

The elliptic genus of a sigma model with Calabi-Yau target is defined to be \cite{Wi}
\begin{equation} \label{ellge}
\phi(\tau,z)=\mathrm{Tr}_{H_{RR}}(q^{L_0-c/24}e( zJ_0)(-1)^F\overline{q}^{\overline{L}_0-c/24})\,,\end{equation}
where $q=e(\tau)$.
Here, $L_0$ resp. $\overline{L}_0$ are the holomorphic resp. antiholomorphic  Virasoro generators defining the grading,  $J_0$ is an operator 
 of $N=2$ supersymmetry defining the charge, and $F$ is
the Fermionic number operator. The central charge $c$ equals $3d$, where $d$ is the
(complex) dimension  of the Calabi-Yau target. $\phi(\tau,z)$ will be independent of $\overline{q}$,
by virtue of the supersymmetry in the Ramond sector, and depends only on the topological
type (rather than the complex structure) of the target.
It is a \textit{weak Jacobi form} of index $m=d/2$ and weight $k=0$ for $\bbZ^2\sdprod
\mathrm{SL}_2(\bbZ)$.
This means that $\phi$ is a holomorphic function of  $z\in\bbC$ and $\tau\in\bbH$, the upper half-plane, which is modular with respect to $\tau$ and quasi-periodic with respect to $z$:
for all $\left({a\atop c}{b\atop d}\right)\in\mathrm{SL}_2(\bbZ)$ and  $l,l'\in\bbZ$,
\begin{eqnarray}
&&\phi\left(\frac{a\tau+b}{c\tau+d},\frac{z}{c\tau+d}\right)=(c\tau+d)^ke( mcz^2/(c\tau+d))
\phi(\tau,z)\ ,\\
&&\phi(\tau,z+l \tau +l')=e(-\I m\,(l^2\tau+2lz))\phi(\tau,z)\ ,\end{eqnarray}
and in addition has Fourier expansion
\begin{equation}
\phi(\tau,z)=\sum_{n\ge 0,l\in\bbZ}c(n,l)q^ny^l\,,\end{equation}
for $y=e(z)$. Consistency of these conditions requires $c(n,l)=(-1)^kc(n,-l)$.
The prefix `{weak}'  denotes that the sum over $n$ is allowed to start at 0, i.e. $\phi$ is merely holomorphic
at the cusp $\tau=\I\infty$; for historical reasons to be a true \textit{Jacobi form} requires a sum over 
$l^2\le 4mn$.

The theory of Jacobi forms is developed in \cite{EZ}.  The definition extends trivially to
weak Jacobi forms 
for other $\bbZ^2\sdprod \Gamma$, where $\Gamma$ is any subgroup of SL$_2(\bbZ)$ of finite 
index. 
Note that because we have insisted here on the full group $\bbZ^2$ of translations, 
$\Gamma$ must be a subgroup of SL$_2(\bbZ)$ because it must act on that
group of translations. This is significant because it means we
will lose the genus-0 property which played so important a role in Monstrous Moonshine
--- recall that many of the Monstrous Moonshine Fuchsian groups aren't subgroups of
SL$_2(\bbZ)$.

Nevertheless, the genus-0 property is what made Monstrous Moonshine special (and rather mysterious), so 
we should seek an analogue for it in Mathieu Moonshine. \cite{CD} have an intriguing proposal:
that the mock modular form $q^{-1/8}\sum_nH_n(g)q^n$ equals a certain regularised 
Rademacher sum. The relation of such a property (in weight-0) to the genus-0 property was
established in \cite{DF}.

The algebra of weak Jacobi forms
for $\bbZ^2\sdprod\Gamma$ of even weight and integral index, is the polynomial
algebra $\mathfrak{m}_\Gamma[\phi_{0,1},\phi_{-2,1}]$, where $\mathfrak{m}_\Gamma$ 
is the ring of holomorphic modular forms for $\Gamma$, and
\begin{equation}
\phi_{0,1}(\tau,z)=4\,\left(\frac{\theta_2(\tau,z)^2}{\theta_2(\tau,0)^2}+\frac{\theta_3(\tau,z)^2}{\theta_3(\tau,0)^2}+ 
\frac{\theta_4(\tau,z)^2}{\theta_4(\tau,0)^2}\right) \,,\ \phi_{-2,1}(\tau,z)=-
\frac{\theta_1(\tau,z)^2}{\eta(\tau)^6} \,,\end{equation}
for the classical Jacobi theta series $\theta_i$ 
 and  Dedekind eta function $\eta$. (For $\Gamma=\mathrm{SL}_2(\bbZ)$ this is Proposition
 9.3 in \cite{EZ}, but the proof for arbitrary finite-index $\Gamma$ is the same.)
 In particular, every weight-$0$ index-1 weak Jacobi form
 $F(\tau,z)$ for $\bbZ^2\sdprod \Gamma$  can be expressed as
 $F=h\phi_{0,1}+f\phi_{-2,1}$, where $h\in\bbC$ and $f$ is a holomorphic weight-$2$ modular
 form for $\Gamma$. For nonlinear $\sigma$-models on K3, this implies the elliptic genus $\phi(\tau,z)$ 
must be proportional to $\phi_{0,1}$ and therefore will be constant throughout the
40-dimensional moduli space of K3 sigma models, something we already knew. 

The worldsheet description is of an $N=2$ vertex operator super algebra (VOSA). The holonomy 
of K3  surfaces extends the $N=2$ to  $N=4$. Therefore,
the space $H_{RR}$  in fact carries a representation of 
the $N=4$ superconformal  algebra at $c=6$, so
\begin{equation}\label{ellgenN4}
\phi_{K3}(\tau,z)=A_{00}ch^s_{1/4,0}(\tau,z)+\sum_{n=0}^\infty A_nch^l_{n+1/4,1/2}(\tau,z)
\,,\end{equation}
where $ch^*_{h,j}(\tau,z)$ is the elliptic genus of the appropriate level 1 $N=4$ representation
($s=\,$short=massless=BPS, $l=\,$long=massive=non-BPS) and $ch^l_{1/4,1/2}$ is to
be interpreted as $ch^s_{0,0}+2ch^s_{1/4,1/2}$.
 Eguchi,  Ooguri \& Tachikawa \cite{EOT} remarked that these coefficients $A_n$, previously
 calculated in \cite{EOTY},
 appear to be dimensions of $M_{24}$-representations (except for $n=0$). More precisely,
 $A_{00}=20$ is the dimension of the virtual representation $H_{00}=
 \rho_1-3\rho_0=\rho_1-3$,  and $A_n$ is the dimension
 of $H_n$, where $H_0$ is the virtual representation $-2$,  and the next few are
{\small\begin{eqnarray}
&&H_1=\rho_2+\overline{\rho_2}\,,\ \  H_2=\rho_3+\overline{\rho_3}\,,\ \  H_3=\rho_7+\overline{\rho_7}\,,\ \  H_4=2\rho_{14}\,,\ \ H_5=2\rho_{19}\,,\nonumber\\
&&H_6=2 \rho_{20}+2\rho_{ 16}\,,\ \ H_7=2\rho_{20}+2\rho_{19}+2\rho_{18}+2\rho_{17}+2\rho_{13}+2\rho_{12}\,,\nonumber\\&&H_8=6 \rho_{20} +2 \rho_{19}+2\rho_{ 18}+4 \rho_{17}+2 \rho_{16}+2 \rho_{15}+2 \rho_{14}
+2\rho_{ 11}+\rho_{10}+\overline{\rho_{10}}+\rho_8+\overline{\rho_8}\,,\nonumber\\
&&H_9=10\rho_{ 20}+8 \rho_{19}+8 \rho_{18}+4 \rho_{17}+4 \rho_{16}+4 \rho_{15}+2 \rho_{14}
+2 \rho_{13}+2 \rho_{12}+ 2\rho_{10}+2\overline{\rho_{10}}\nonumber\\
&&\qquad\qquad +2 \rho_9+2 \rho_7+2\overline{\rho_7}+ 2 \rho_6\,.\nonumber
\end{eqnarray}}
(Incidentally, \cite{EHb} note that expanding the elliptic genus into $N=2$ characters (rather
than $N=4$) seems to carry a representation of PSL$_2(\bbZ_{11})$, a simple group of
order 660.)

This observation of \cite{EOT} is very reminiscent of the Monstrous Moonshine observation of McKay 
who noted that the expansion coefficients of the Hauptmodul $J$-function seem to be dimensions of Monster 
group $\bbM$ representations $V_n$ (see e.g. \cite{Ga} for a review and references). 
Thompson suggested considering what are now called the \textit{McKay-Thompson series}, defined
formally  for all $g\in\bbM$ by
\begin{equation} T_g(\tau)=\sum_{n=-1}^\infty \mathrm{ch}_{Vn}(g)\,q^n\,,\end{equation}
where $q=e^{2\pi\I\tau}$. Note that these $T_g$ are  constant on each of the 194 conjugacy
classes of $\bbM$ --- they are called \textit{class functions}.
After staring at their first few coefficients, Conway \& Norton \cite{CN} conjectured that
these $T_g$ are modular functions (in fact Hauptmoduls=normalised generators of the field of modular functions)
for some genus-0 Fuchsian group $\Gamma_g$ commensurable with (but not in general a
subgroup of) SL$_2(\bbZ)$.  To make this conjecture more precise, we should turn the logic around and associate to each conjugacy
class $K_g$ in $\bbM$ a uniquely specified modular function $T_g$. From their Fourier
coefficients  we obtain 
a sequence $\chi_{n}$ of \textit{class functions}, i.e. linear combinations over $\bbC$
of irreducible characters of $\bbM$. Because the coefficients of the
$T_g$ are rational (in fact integral), it is immediate that this linear combination can be
taken over $\bbQ$,  but it is very hard to show they can be taken over $\bbZ$ (i.e. that these
class functions $\chi_n$ are {virtual characters}) and even harder to show that they are
combinations  over $\bbZ_{\ge 0}$ (i.e. that
the $\chi_n$ are actually characters of $\bbM$-representations $M_n$). 
Atkin, Fong \& Smith (see \cite{Sm}) tried to prove with the aid of a computer 
that these $\chi_n$ are indeed true characters, and managed to reduce the proof
to a fairly plausible statement they called Conjecture 2.3 (it is often claimed in
the literature --- see e.g. \cite{Ga} --- that \cite{Sm} proves the $\chi_n$ are true
characters, but this is false). That Conway-Norton conjecture was finally proved in
\cite{Bo}, independently of \cite{Sm}; thanks to this work,  
 we now know much more: there is a vertex 
operator algebra (VOA) $\cV^\natural$ (constructed in \cite{FLM}) with automorphism group $\bbM$,  
whose \textit{twining characters} $\mathrm{Tr}_{\cV^\natural}gq^{L_0-1}$ equal the McKay--Thompson
series $T_g$ (the `$-1$' in the exponent is the usual shift by $c/24$). In other words,  the $\bbM$-modules $V_n$ are the eigenspaces of
$V^\natural$ with respect to the Virasoro operator $L_0$. (The adjective `twining' is short
for `intertwining', and is used as an alternative to the over-used word `twisted'.)

Thus we are led to test further
this Mathieu Moonshine observation of \cite{EOT}, by introducing for each $g\in M_{24}$
\begin{equation}
\phi_g(\tau,z)=\mathrm{ch}_{H00}(g)\,\mathrm{ch}_{1/4,0}(\tau,z)+\sum_{n=0}^\infty \mathrm{ch}_{Hn}(g)\,\mathrm{ch}_{n+1/4,1/2}(\tau,z)
\,.\end{equation}
These should be the \textit{twining elliptic genera} 
\begin{equation}
\phi_g(\tau,z)=\mathrm{Tr}_{H_{RR}}gq^{L_0-c/24}e( zJ_0)(-1)^F\overline{q}^{\overline{L}_0-c/24}
\,,\label{tweg}\end{equation}
although that expression is purely formal as none of these K3 sigma models will have a
symmetry consisting of all of $M_{24}$. In fact \cite{Mu}, the group of symplectic automorphisms of
any K3 surface will typically be trivial, will never exceed order 960, and will only contain
elements of order $\le 8$.  And no automorphism, symplectic or otherwise, of a K3 surface
can have order 23. But ignore these subtleties for now.
Twining elliptic genera should have good modular properties --- we would expect them to be
weight-0 index-1 Jacobi forms for  $\bbZ^2\sdprod\Gamma_0(|g|)$ (recall we write $|g|$
for the order of $g$) though possibly with multiplier. This is enough to guess 
 with effort these $\phi_g$ from
the first few coefficients (just as was done in Monstrous Moonshine by Conway \& Norton). 

Conjectured expressions for $\phi_g$, for all $g\in M_{24}$ (constant on conjugacy classes of course) 
 were  obtained in \cite{Ch, GHV1,GHV2,EH}. We have 
\begin{equation}\label{phigfg}
\phi_g(\tau,z)=\frac{w_g}{12}\phi_{0,1}(\tau,z)+f_g(\tau)\,\phi_{-2,1}(\tau,z)\,,
\end{equation}
{where $w_g\in\bbZ$ is the \textit{Witten index} given in Table 2, and $f_g(\tau)$ is some holomorphic
weight-2 modular form (with trivial multiplier) for $\Gamma_0(|g|h_g)$ for $h_g$ in Table 2.}
$w_g$ equals the character value of $\rho_1+1$
at $g$; note that it vanishes iff $K_g$ doesn't intersect $M_{23}$ --- this isn't  deep: $\rho_1+1$ 
is the permutation representation, so $K_g$ intersects $M_{23}$ iff $\rho_1+1$ has a fixed point,
iff its character value (which equals the number of fixed-points) is nonzero. $h_g$ is the length of 
the shortest cycle in the cycle shape of $g$.
Note that in all cases, $h_g|\mathrm{gcd}(N,12)$ (where $N$ is the order of the element), and $h_g=1$ iff the conjugacy class $K_g$
intersects $M_{23}$.  We include for later use the order of the centraliser $C_{M24}(g)$
of $g$ in $M_{24}$, and the index of $\Gamma_0(|g|)$ in SL$_2(\bbZ)$ (the index of
$\Gamma_0(\prod_ip_i^{a_i})$, where $p_i$ are distinct
primes,  is $\prod_i(p_i+1)p_i^{a_i-1}$), as well as the bound $a\times 10^b$ obtained
by Theorem 3 below.
In Table 2 and elsewhere, we write $7AB$  for $7A\cup 7B$ etc; we collect together conjugacy classes like this 
when the corresponding Jacobi forms $\phi_g$ coincide. 

 \medskip
\centerline {{\bf Table 2.} Data for $M_{24}$ conjugacy classes}
$${\scriptsize\vbox{\tabskip=0pt\offinterlineskip
  \def\tablerule{\noalign{\hrule}}
  \halign to 6.9in{
    \strut#&
    \hfil#&\vrule#&\vrule#&\hfil#&\vrule#&    
\hfil#&\vrule#&\hfil#&\vrule#&\hfil#&\vrule#&\hfil#&\vrule#&\hfil#&\vrule#&\hfil#&\vrule#&\hfil#&\vrule#&\hfil#&\vrule#&\hfil#&\vrule#&\hfil#&\vrule#&\hfil#&\vrule#&\hfil#&\vrule#&\hfil#&\vrule#&\hfil#&\vrule#&\hfil#&\vrule#&\hfil#&\vrule#&\hfil#&\vrule#&\hfil#&\vrule#&\hfil#&\vrule#&\hfil#
\tabskip=0pt\cr
&$K_g\,\,$&\,&&\,\,\,\,\,1A\,\,\,\,&&\,\,\,\,\,2A\,\,\,\,\,&&\,\,\,\,2B\,\,\,\,&&\,\,\,3A\,\,\,&&\,\,3B\,\,&&\,\,4A\,\,&&\,\,4B\,\,&&\,\,4C\,\,&&\,\,5A\,\,&&\,\,6A\,\,&&\,\,6B\,\,&&\,7AB\,&&\,\,8A\,\,&&\,\,10A\,\,&&\,\,11A\,\,&&\,\,12A\,\,&&\,\,12B\,\,&&\,14AB\,&&\,15AB\,&&\,21AB\,&&\,23AB\cr
\tablerule&$h_g\,\,$&\,&&\,\,1\,\,&&\,\,1\,\,&&\,\,2\,\,&&\,\,1\,\,&&\,\,3\,\,&&\,\,2\,\,&&\,\,1\,\,&&\,\,4\,\,&&\,\,1\,\,&&\,\,1\,\,&&\,\,6\,\,&&\,\,1\,\,&&\,\,1\,\,&&\,\,2\,\,&&\,\,1\,\,&&\,\,2\,\,&&\,\,12\,\,&&\,\,1\,\,&&\,\,1\,\,&&\,\,3\,\,&&\,\,1\,\,\,\cr
\tablerule&$w_g\,\,$&\,&&\,\,24\,\,&&\,\,8\,\,&&\,\,0\,\,&&\,\,6\,\,&&\,\,0\,\,&&\,\,0\,\,&&\,\,4\,\,&&\,\,0\,\,&&\,\,4\,\,&&\,\,2\,\,&&\,\,0\,\,&&\,\,3\,\,&&\,\,2\,\,&&\,\,0\,\,&&\,\,2\,\,&&\,\,0\,\,&&\,\,0\,\,&&\,\,1\,\,&&\,\,1\,\,&&\,\,0\,\,&&\,\,1\,\,\,\cr
\tablerule&$|C(g)|$&\,&&$|M_{24}|\,$&&21504\,&&7680\,&&1080\,&&504\,&&384\,&&128\,&&96\,&&60\,&&24\,&&24\,&&42\,&&16\,&&20\,&&11\,&&12\,&&12\,&&14\,&&15\,&&21\,&&23\,\,\,\cr
\tablerule&index\,&\,&&\,\,1\,\,&&\,\,3\,\,&&\,\,3\,\,&&\,\,4\,\,&&\,\,4\,\,&&\,\,6\,\,&&\,\,6\,\,&&\,\,6\,\,&&\,\,6\,\,&&\,\,12\,\,&&\,\,12\,\,&&\,\,8\,\,&&\,\,12\,\,&&\,\,18\,\,&&\,\,12\,\,&&\,\,24\,\,&&\,\,24\,\,&&\,\,24\,\,&&\,\,24\,\,&&\,\,32\,\,&&\,\,24\,\,\,\cr
\tablerule&$a\!\times\! 10^b$&\,&&3\,\,6\,\,&&5.7\,6\,\,&&2.1\,9\,&&4.6\,7\,&&\,2.3\,12&&\,2.9\,11&&\,2.1\,8\,&&5.1\,14&&\,6.8\,8\,&&\,2.9\,9\,&&\,1.5\,18\,&&7.5\,9\,&&\,2.2\,10\,&&\,3.5\,14\,&&\,2.6\,11&&\,1.5\,15\,&&\,5.8\,23\,&&1.7\,12\,&&2.9\,12\,&&8.7\,18\,&&\,8.2\,13\,\cr
\noalign{\smallskip}}}}$$

Extracting $f_g$ from  \eqref{phigfg}, we obtain
\begin{eqnarray}
&&f_g(\tau)=w_g\,\frac{\theta_3^4
+\theta_4^4}{12}-w_g\,{\theta_3\theta_4}\left(\frac{1}{4}+\sum_{n=1}^\infty\frac{q^{n(n+1)/2}}{1+
q^n}\right)-\eta^3q^{-1/8}\sum_{n=0}^\infty H_n(g)q^n\nonumber\\
&&\qquad\quad=-\eta^3q^{-1/8}\sum_{n=0}^\infty\left(H_n(g)-\frac{w_g}{24}H_n(1)\right)q^n\,.\label{eq:fg}\end{eqnarray}
These modular forms $f_g(\tau)$ 
are given explicitly (i.e. independently of knowing the $H_n(g)$) for all $g\in M_{24}$ in \cite{GHV3,EH} (see also Table 2 of \cite{CD}). Indeed, recall the Eisenstein series $E_2(\tau)
=1-24q-72q^2-\cdots$, a quasi-modular form for SL$_2(\bbZ)$. Then for each integer $n>1$, 
\begin{equation} \label{E2n}
E_2^{(n)}(\tau):=\frac{1}{n-1}(E_2(n\tau)-E_2(\tau))=1+\frac{24}{n-1}\sum_{k=1}^\infty
\sigma_1(k)\left(q^k-q^{nk}\right)\end{equation}
 is a holomorphic modular form of weight 2
for $\Gamma_0(n)$ with trivial multiplier ($\sigma_1(k)=\sum_{d|k}d$). Writing $\eta(n):=\eta(n\tau)$ and $\eta=\eta(\tau)$,
we have: $f_{1A}=0$,
{\small\begin{eqnarray} 
f_{2A}&=&\frac{4}{3}E_2^{(2)}=\frac{4}{3}+32q+32q^2+128q^3+32q^4+192q^5+128q^6 +256q^7+32q^8+416q^9\cdots   \,,\nonumber\\
f_{2B}&=&2\frac{\eta^8}{\eta(2)^4}= 2-16q+48q^2-64q^3+48q^4-96q^5+192q^6 -128q^7 +48q^8- 208 q^9 +\cdots \,,\nonumber\\
f_{3A}&=&\frac{3}{2}E_2^{(3)}=\frac{3}{2}+18q+54q^2+18q^3+126q^4+108q^5+54q^6 +144q^7+270q^8 + 18 q^9+\cdots  \,,\nonumber\\
f_{3B}&=& 2\frac{\eta^6}{\eta(3)^2}=2 - 12 q  + 18q^2 + 24q^3  - 84q^4  + 36q^5  + 72q^6-96q^7+90q^8+ 24 q^9+\cdots  \,,\nonumber\\
f_{4A}&=&2\frac{\eta(2)^8}{\eta(4)^4}=2 - 16q^2 + 48q^4 - 64q^6  +48q^8 +0q^9+\cdots \,,\nonumber\\
f_{4B}&=&  2E_2^{(4)}-\frac{1}{3}E_2^{(2)}=\frac{5}{3} + 8q+ 40q^2+ 32q^3+ 40q^4+ 48q^5+ 160q^6 +64q^7+40q^8+ 104 q^9+\cdots \,,\nonumber\\
f_{4C}&=& 2\frac{\eta^4\eta(2)^2}{\eta(4)^2}=2 - 8q + 32q^3- 16 q^4- 48q^5 +64q^7 +48q^8- 104 q^9+\cdots \,,\nonumber
\end{eqnarray} }{\small\begin{eqnarray} 
f_{5A}&=& \frac{5}{3}E_2^{(5)}=\frac{5}{3} + 10q+ 30q^2+ 40q^3+ 70q^4+ 10q^5+ 120 q^6+80q^7+150q^8 + 130 q^9 +\cdots \,,\nonumber\\
f_{6A}&=&\frac{5}{2}E_2^{(6)}-\frac{1}{2}E_2^{(3)}-\frac{1}{6}E_2^{(2)}=   \frac{11}{6 }+ 2q+ 14q^2+ 26q^3+ 38q^4+ 12q^5+ 38 q^6+16q^7+86q^8 + 98 q^9+\cdots \,,\nonumber\\
f_{6B}&=& 2\frac{\eta^2\eta(2)^2\eta(3)^2}{\eta(6)^2}=2 - 4q- 6q^2+ 8q^3+ 12q^4+ 12q^5- 24q^6 -32q^7-6q^8 + 8 q^9+\cdots \,,\nonumber\\
f_{7AB}&=& \frac{7}{4}E_2^{(7)}=\frac{7}{4} + 7q+ 21q^2+ 28q^3+ 49q^4+ 42q^5 + 84q^6+7q^7+105q^8 + 91 q^9  +\cdots \,,\nonumber\\
f_{8A}&=& \frac{7}{3}E_2^{(8)}-\frac{1}{2}E_2^{(4)}=\frac{ 11}{6} + 4q+ 12q^2+ 16q^3+ 44q^4+ 24q^5+ 48q^6+32q^7+44q^8+ 52 q^9 +\cdots \,,\nonumber\\
f_{10A}&=&  2\frac{\eta^3\eta(2)\eta(5)}{\eta(10)}=2 - 6q- 2q^2+ 16q^3- 2q^4- 6 q^5- 8q^6-8q^7-2q^8+ 2 q^9 +\cdots \,,\nonumber\\
f_{11A}&=&\frac{11}{6}E_2^{(11)}-\frac{22}{5}\eta^2\eta(11)^2=\frac{11}{6} + 22q^2+ 22q^3+ 22q^4+ 22q^5+ 44 q^6 +44q^7+66q^8+ 66 q^9  +\cdots \,,\nonumber\\
f_{12A}&=& 2\frac{\eta^3\eta(4)^2\eta(6)^3}{\eta(2)\eta(3)\eta(12)^2}=2 - 6q+ 2q^2+ 6q^3- 6q^4 + 12q^5 - 10 q^6-6q^8 - 6 q^9 +\cdots \,,\nonumber\\
f_{12B}&=& 2\frac{\eta^4\eta(4)\eta(6)}{\eta(2)\eta(12)}=  2 - 8q+ 6q^2+ 8q^3- 4q^4- 24q^6+16q^7-6q^8+ 16 q^9+\cdots \,,\nonumber\\
f_{14AB}&=&  \frac{91}{36}E_2^{(14)}-\frac{7}{12}E_2^{(7)}-\frac{1}{36}E_2^{(2)}-\frac{14}{3}\eta\eta(2)\eta(7)\eta(14)\nonumber\\&=&\frac{23}{12}- 3q+ 11q^2+ 16q^3+ 11q^4+ 10 q^5+ 16q^6+25q^7+39q^8 + 17 q^9+\cdots \,,\nonumber \\
f_{15AB}&=&  \frac{35}{16}E_2^{(15)}-\frac{5}{24}E_2^{(5)}-\frac{1}{16}E_2^{(3)}-\frac{15}{4}\eta\eta(3)\eta(5)\eta(15)\nonumber\\&=& \frac{23}{12} - 2q+ 9q^2+ 13q^3 + 16q^4+ 13q^5+ 24q^6+ 14q^7+15q^8+ 28 q^9+\cdots \,,\nonumber\\
f_{21AB}&=&  \frac{7}{3}\frac{\eta^3\eta(7)^3}{\eta(3)\eta(21)}-\frac{1}{3}\frac{\eta^6}{\eta(3)^2}
=2 - 5q- 3q^2+ 10q^3 + 7q^4- 6q^5- 12q^6- 5 q^7+6q^8 + 10 q^9+\cdots \,,\nonumber\\
f_{23AB}&=&  \frac{23}{12}E_2^{(23)}-69\eta^2\eta(23)^2+
92\eta(2)^2\eta(46)^2+92\eta\eta(2)\eta(23)\eta(46)+23\frac{\eta^3\eta(23)^3}{\eta(2)\eta(46)}
\nonumber\\&=& \frac{23}{12}+ 23q^3+ 23q^4 + 23q^6 + 46q^8+ 23 q^9+\cdots \,.\nonumber\end{eqnarray} } 

Equation (7.16) of \cite{DMZ} tells us
\begin{equation}\label{dmz}
\eta^3q^{-1/8}\sum_{n=0}^\infty H_n(1)\,q^n=48F^{(2)}(\tau)-2E_2(\tau)\,,\end{equation}
where $F^{(2)}_2(\tau)=\sum_{n>m>0,n\not\equiv_2m}(-1)^nmq^{mn/2}$.
We read off from \eqref{E2n} and \eqref{dmz} that
\begin{align}\label{E22}
E^{(2)}_2&\equiv 1\ \ (\mathrm{mod}\ 24)\,,\\
E^{(3)}_2&\equiv 1\ \ (\mathrm{mod}\ 12)\,,\label{E23}\\
\label{3.2}
\eta^3q^{-1/8}\sum_{n=0}^\infty H_n(1)\,q^n&\equiv -2\qquad (\mathrm{mod}\ 48)\,.
\end{align}
These will be used in the Section 3 proofs. (We thank Thomas Creutzig and Gerald H\"ohn
for bringing \eqref{dmz} and its consequence \eqref{3.2} to our attention.)

 \medskip\noindent\textbf{Lemma 1.} \cite{GHV2} \textit{Let $g\in M_{24}$ have order $|g|$, with
 parameter $h$ given in Table 2. Then the twining character $\phi_g(\tau,z)$ is a Jacobi
 form of index 1 and weight 0 under $\Gamma_0(|g|)$ with multiplier defined by}
\begin{equation}
\phi_g\left(\frac{a\tau+b}{c\tau+d},\frac{z}{c\tau+d}\right) =e^{2\pi\I cd/(|g|\,h)}e^{2\pi\I cz^2/(c\tau+d)}
\phi_g(\tau,z)\,.
\label{phin}\end{equation}
\textit{The function $f_g(\tau)$is a holomorphic modular form of weight 2  for $\Gamma_0(|g|)$
 with multiplier $e^{2\pi\I cd/(|g|\,h)}$, i.e. $f_g\left(\frac{a\tau+b}{c\tau+d}\right)=(c\tau+d)^2e^{2\pi\I cd/(|g|\,h)}f_g(\tau)$ $\forall \left({a\atop c}{b\atop d}\right)\in\Gamma_0(|g|)$.}\medskip

$\Gamma_0(n)$ as usual consists of all matrices $\left(\begin{matrix} {a\atop c}{b\atop d}\end{matrix}\right)\in\mathrm{SL}_2(\bbZ)$ with $n|c$. 
 The appearance of both $\Gamma_0(|g|)$ and the multiplier in Lemma 1 are
not at all mysterious, as explained in the concluding section.

 Again, reversing the logic, we can use these
 conjectured expressions for $\phi_g$ to define class functions $H_{00}(g)$ and $H_n(g)$.
 It 
was checked by brute force that these class functions are in fact true characters, for
all $n\le 500$ (except $H_0$ and $H_{00}$ which are merely virtual).
This is our Theorem A.

 One would certainly want more than weak Moonshine (Theorem A)
--- e.g. we should have an $M_{24}$-worth of twisted twining elliptic genera with nice modularity.
This is now done \cite{Ron}.
Unlike Monstrous Moonshine, where the algebraic content comes from the VOA $V^\natural$,
the algebraic (not to mention geometric and physical) meaning of K3-Mathieu Moonshine is still unclear.  
Although we have learned much in the two years or so since \cite{EOT}, we still don't really
know the right questions to ask.

\section{Weak  K3-Mathieu Moonshine I: Integrality}

In this section and the next, we prove Theorems A and B, which were stated in
Section 2. Our proof falls into two independent steps: \textit{integrality} (this section)
and \textit{positivity} (next section).

Let $G$ be any finite group. Let $\cP_G$ be the set of all pairs $(p,K_g)$ where $p$ is a prime
dividing the order $|G|$, $K_g$ is a \textit{$p$-regular} conjugacy class in $G$ (i.e. the order
$|\langle g\rangle|$ of $g$ is coprime to $p$), and $p$ divides the order $|C_G(g)|$ of
the centraliser. There are precisely 22 pairs in $\cP_G$ for $G=M_{24}$, which we list in Table 3.
For reasons to become clear shortly, we include the highest power $p^\pi$ of $p$ dividing
$|C_G(g)|$, and the \textit{$p'$-section} $S$ which we define next paragraph.

 \medskip
\centerline {{\bf Table 3.} Data for virtual character proof}
$${\scriptsize\vbox{\tabskip=0pt\offinterlineskip
  \def\tablerule{\noalign{\hrule}}
  \halign to 6.5in{
    \strut#&
    \hfil#&\vrule#&\vrule#&\hfil#&\vrule#&    
\hfil#&\vrule#&\hfil#&\vrule#&\hfil#&\vrule#&\hfil#&\vrule#&\hfil#&\vrule#&\hfil#&\vrule#&\hfil#&\vrule#&\hfil#&\vrule#&\hfil#&\vrule#&\hfil#&\vrule#&\hfil#&\vrule#&\hfil#&\vrule#&\hfil#&\vrule#&\hfil#&\vrule#&\hfil#&\vrule#&\hfil#&\vrule#&\hfil#&\vrule#&\hfil#&\vrule#&\hfil#&\vrule#&\hfil#&\vrule#&\hfil#
\tabskip=0pt\cr
&$p\,\,$&\,&&\,\,2\,\,&&\,\,2\,\,&&\,\,2\,\,&&\,\,2\,\,&&\,\,2\,\,&&\,\,2\,\,&&\,\,3\,\,&&\,\,3\,\,&&\,\,3\,\,&&\,\,3\,\,&&\,\,3\,\,&&\,\,3\,\,&&\,\,3\,\,&&\,\,3\,\,&&\,\,5\,\,&&\,\,5\,\,&&\,\,5\,\,&&\,\,7\,\,&&\,\,7\,\,&&\,\,7\,\,&&\,\,11\,\,&&\,\,23\cr\tablerule 
&$K_g\,\,$&\,&&\,\,1A\,\,&&\,\,3A\,\,&&\,\,3B\,\,&&\,\,5A\,\,&&\,\,7A\,\,&&\,\,7B\,\,&&\,\,1A\,\,&&\,\,2A\,\,&&\,\,2B\,\,&&\,\,4A\,\,&&\,\,4C\,\,&&\,\,5A\,\,&&\,\,7A\,\,&&\,\,7B\,\,&&\,\,1A\,\,&&\,\,2B\,\,&&\,\,3A\,\,&&\,\,1A\,\,&&\,\,2A\,\,&&\,\,3B\,\,&&\,\,1A\,\,&&\,\,1A\cr\tablerule \tablerule 
&$\pi\,\,$&\,&&\,\,10\,\,&&\,\,3\,\,&&\,\,3\,\,&&\,\,2\,\,&&\,\,1\,\,&&\,\,1\,\,&&\,\,3\,\,&&\,\,1\,\,&&\,\,1\,\,&&\,\,1\,\,&&\,\,1\,\,&&\,\,1\,\,&&\,\,1\,\,&&\,\,1\,\,&&\,\,1\,\,&&\,\,1\,\,&&\,\,1\,\,&&\,\,1\,\,&&\,\,1\,\,&&\,\,1\,\,&&\,\,1\,\,&&\,\,1\cr \tablerule 
&$S\,\,$&\,&&\,\,1A,2AB\,\,&&\,\,3A\,\,&&\,\,3B\,\,&&\,\,5A\,\,&&\,\,7A\,\,&&\,\,7B\,\,&&\,\,1A\,\,&&\,\,2A\,\,&&\,\,2B\,\,&&\,\,4A\,\,&&\,\,4C\,\,&&\,\,5A\,\,&&\,\,7A\,\,&&\,\,7B\,\,&&\,\,1A\,\,&&\,\,2B\,\,&&\,\,3A\,\,&&\,\,1A\,\,&&\,\,2A\,\,&&\,\,3B\,\,&&\,\,1A\,\,&&\,\,1A\cr  
&&\,&&\,\,4ABC\,\,&&\,\,6A\,\,&&\,\,6B\,\,&&\,\,10A\,\,&&\,\,14A\,\,&&\,\,14B\,\,&&\,\,3A\,\,&&\,\,6A\,\,&&\,\,6B\,\,&&\,\,12A\,\,&&\,\,12B\,\,&&\,\,15A\,\,&&\,\,21A\,\,&&\,\,21B\,\,&&\,\,5A\,\,&&\,\,10A\,\,&&\,\,15A\,\,&&\,\,7A\,\,&&\,\,14A\,\,&&\,\,21A\,\,&&\,\,11A\,\,&&\,\,23A\cr  
&&\,&&\,\,8A\,\,&&\,\,12A\,\,&&\,\,12B\,\,&&&&&&&&\,\,3B\,\,&&&&&&&&&&\,\,15B\,\,&&&&&&&&&&\,\,15B\,\,&&\,\,7B\,\,&&\,\,14B\,\,&&\,\,21B\,\,&&&&\,\,23B\cr \tablerule
&$q^n\,\,$&\,&&\,\,$q^4$\,\,&&\,\,$q^8$\,\,&&\,\,$q^{48}$\,\,&&\,\,$q^6$\,\,&&\,\,$q^4$\,\,&&\,\,
$q^4$\,\,&&\,\,$q^2$\,\,&&\,\,$q^2$\,\,&&\,\,$q^2$\,\,&&\,\,$q^4$\,\,&&\,\,$q^6$\,\,&&\,\,$q^4$\,\,&&
\,\,$q^{16}$\,\,&&\,\,$q^{16}$\,\,&&\,\,$q^1$\,\,&&\,\,$q^3$\,\,&&\,\,$q^4$\,\,&&\,\,$q^1$\,\,&&\,\,$q^4$\,
\,&&\,\,$q^5$\,\,&&\,\,$q^2$\,\,&&\,\,$q^4$\cr
 \noalign{\smallskip}}}}$$

For any $h\in G$ of order $n=p^km$, where gcd$(p,m)=1$, find $a,b\in\bbZ$ such that
$1=ap^k+bm$. Then $h=h_{p'}h_p$ where the \textit{$p'$-part} $h_{p'}=h^{ap^k}$ has order $m$
and the \textit{$p$-part} $h_p=h^{bm}$ has order $p^k$. The $p'$-section $S$ for the pair
$(p,K_g)\in\cP_G$ is defined  to be $S=\{h\in G\,|\, h_{p'}\in K_g\}$. $S$ is clearly the union of
conjugacy classes; in Table 3 we list those conjugacy classes.

Suppose we have an integer-valued class function $c:G\rightarrow\bbZ$, i.e. $c$ is constant
on conjugacy classes. We are interested in  $H$ being one of the head characters $H_n$, and we
want to prove it is the character of an $M_{24}$-representation. The hard part is to prove that
$A$ is a virtual character, i.e. a linear combination over $\bbZ$ of irreps. The starting point is the
following characterisation of virtual characters, attributed to Thompson and based on the classical
theorem of Brauer. It is a refinement of the following basic idea, which the reader can
verify for himself: If $\chi$ is an integer-valued character and $g$ is an element of order a prime 
power $p^\pi$, then $\chi(g)\equiv \chi(e)$ (mod $p$).

Choose any $(p,K_g)\in \cP_G$ and write the $p'$-section $S$ as a disjoint union $\cup_{i=1}^kK_i$ of conjugacy classes, for some $k$ depending on $(p,S)$. 
Let $R_p$ be the tensor product $\widehat{\bbZ}_p\otimes_\bbZ \bbZ[\xi_{|G|}]$, where
$\widehat{\bbZ}_p$ are the $p$-adic integers and $\xi_n=e^{2\pi\I/n}$. $R_p$ is introduced 
to make formal sense of the orthogonality relations we're about to introduce;
the ring  $\bbZ[\xi_{|G|}]$ arises because all irreducible $G$-characters take values there.
Define $\cM_{(p,K_g)}$ to be the set of all $k$-tuples $(\ell_1,\ldots,\ell_k)\in R_p^k$ such that
$\sum_{i=1}^k\ell_i\chi(K_i)\in p^\pi R_p$ for all irreducible characters $\chi$. If $c$ is a virtual character, then $c$ is likewise orthogonal to
$\cM_{(p,S)}$. What is important to us is the converse:

\medskip\noindent\textbf{Lemma 2.} \cite{Sm} \textit{Let $R_p$ and $\cM_{(p,K_g)}$ be as above.
Let $c:G\rightarrow \bbZ$ be an integer-valued class function of $G$. Then $c$ is a virtual
character of $G$, i.e. 
a linear combination over $\bbZ$ of irreducibles, if $\sum_{i=1}^k \ell_ic(K_i)\in
p^\pi R_p$ for all $(\ell_1,\ldots,\ell_k)\in\cM_{(p,K_g)}$ and all $(p,K_g)\in\cP_G$.}\medskip

Of course we want to apply this to $c=H_n$, for each $n>0$. 
The reason this characterisation
is helpful is that much is known about reductions of modular forms modulo powers of primes,
as we will see.

Suppose $f(\tau)=\sum_{n=0}^\infty f_n q^n,g(\tau)=\sum_{n=0}^\infty g_nq^n$ are holomorphic
modular forms of weight $k\in\bbZ$ for some finite-index subgroup $\Gamma$ and 
 multiplier $\mu:\Gamma
\rightarrow \bbC^\times$, and $\mu$ has finite
order $M$. If the Fourier coefficients $f_n,g_n$ are equal for all $n\le k\|\mathrm{SL}_2(\bbZ)/\Gamma\|
/12$, then $f=g$. This is an immediate consequence of the classical valence formula for
$\Gamma$ applied to $(f-g)^M$. Incidentally, because of this and Lemma 1, the $q$-expansions
provided earlier are more than enough to uniquely specify $f_g$
 (and hence $\phi_g$)  --- in fact the coefficients
$n\le 5$ suffice. What is much more surprising is that this also applies to the modular
forms mod prime powers:

\medskip\noindent\textbf{Lemma 3.} \textit{Suppose $f(\tau)=\sum_{l\ge 0} f_lq^l\in\bbQ[[q^{1/N}]]$ 
and $g(\tau)
=\sum_{l\ge 0} g_lq^l\in\bbQ[[q^{1/N}]]$ are holomorphic modular forms of rational weight $k\in\bbQ$ and
with some multiplier $\mu$ for some subgroup $\Gamma$
of SL$_2(\bbZ)$. Let $m:=\|\mathrm{SL}_2(\bbZ)/\Gamma\|$ be the index. We require  that $\mu$ has finite order,
i.e. that all values $\mu(\gamma)$ for $\gamma\in\Gamma$ are $M$th roots of 1 for some
$M$.} 

\smallskip\noindent\textbf{(a)} \textit{Suppose the Fourier coefficients $f_l$ and $g_l$ are integral for all $l$.
Suppose we have an integer $n>0$ such that  $f_l\equiv g_l$ (mod $n$) for all $l\le km/12$.
Then $f_l\equiv g_l$ (mod $p$) for all $l$.}

\smallskip\noindent\textbf{(b)} \textit{Suppose now that the weight $k$ is integral and that the kernel of $\mu$ 
is a congruence  subgroup, i.e. contains $\Gamma(N')$ for
some $N'$. If the coefficients $f_l$ are integral for all $0\le l\le km/12$,
then $f_l\in\bbZ$ for all $l$.}
\medskip

\noindent\textit{Proof.} The key result upon which this Lemma rests is Sturm's Theorem \cite{St},
which is part (a) in the special case where the multiplier $\mu$ is trivial and the weight is integral
and $n$ is a prime (Sturm however allows the Fourier coefficients to be algebraic integers).  

It suffices to prove the Lemma for $g=0$. Start with part (a). Let $p$ be any prime dividing $n$. Choose
an integer $K>0$ so that $Kk\in\bbZ$ and $\mu^K$ is identically 1. Then $F:=f^K$ is a 
holomorphic modular form for $\Gamma$ of integer weight $Kk$ and trivial multiplier,
whose Fourier coefficients $F_l$ are multiples of the prime $p$ for all $l\le Kkm/12$ (since
each $F_l$ will be a sum of monomials in $f_{l_i}$'s at least one of which has $l_i\le km/12$).
Then Sturm's theorem applies, and we find that $p$ divides \textit{all} coefficients $F_l$.
This must imply that $p$ divides all coefficients of $f$ (otherwise, choose the smallest $l$
such that $p$ doesn't divide $f_l$, and note that $p$ would fail to divide $F_{Kl}$). This
means that $f/p$ will have integer coefficients, so $f/p$ will obey all hypotheses of part (a) with now
$n$ replaced with $n/p$. Repeat with another prime dividing $n/p$, until you've reduced
$n$ to 1.

 To see Lemma 3(b), let $V$ be the space of all holomorphic weight $k$ modular forms
 for $\Gamma(N')$ with trivial multiplier. $V$ is finite-dimensional with an
 integral basis $f^{(i)}\in\bbZ[[q^{1/N'}]]$.
 Therefore we can write $f=\sum_i a_if^{(i)}$ where $a_i\in\bbQ$. This means the coefficients
 of $f$ have denominator bounded by the lcm of all denominators of the $a_i$. Let $L$
 be the smallest positive integer such that $Lf\in\bbZ[[q]]$. If $p$ is a prime dividing $L$ then
 $p$ will divide the first $mk/12$ coefficients of $Lf$ by the hypothesis on $f$, hence will divide all coefficients by
Lemma 3(a), contradicting the minimality of $L$. This means $L=1$, and we're done.\qquad\textbf{QED to Lemma 3}

\medskip\noindent\textbf{Theorem 1.} \textit{Each $H_n$ ($n\ge 1$) is a virtual character of $M_{24}$
(i.e. a linear combination over $\bbZ$ of irreducible characters of $M_{24}$).}\medskip

\noindent\textit{Proof.}  First, we need to show the class functions $H_n$ are integer-valued. 
We will do this by studying the $f_g$'s.
One complication is that the coefficients of $f_g(\tau)$ are not quite integers, as can be seen
by the displayed values in Section 2. The problem is the constant term $\frac{1}{4}$ in
the middle term  of \eqref{eq:fg}, and the
12 in the denominator of its first term. In fact only the constant term of $w_g\,(\theta_3^4
+\theta_4^4)/{12}-w_g\theta_3\theta_4/{4}$ can
fail to be integral. To see this, it suffices to show that $\theta_3\theta_4\in 1+4q\bbZ[[q]]$
and $\theta_3^4+\theta_4^4\in 2+12q\bbZ[[q]]$. That these both hold modulo 4 follows
from $\theta_3=1+2\sum+2\sum'$ and $\theta_4=1+2\sum-2\sum'$, using obvious notation.
That $\theta_3^4+\theta_4^4\equiv 1$ (modulo 3) follows from  comparing $(\theta_3^4+
\theta_4^4)^2$ and $E_4$ modulo 3, using Lemma 3 with $\Gamma=\Gamma(2)$. 

From Lemma 3(b) and checking the integrality of the first few coefficients of $f_g$ we learn now 
that all coefficients of $f_g-w_g\,(\theta_3^4
+\theta_4^4)/{12}+w_g\theta_3^2\theta_4^2/{4}$ are integral (the theta function contribution
kills the fractional part of the constant term). This means that all coefficients of each $f_g$ are
integral except possibly for the constant term, and also (from the first equality in \eqref{eq:fg})
that all coefficients of
$\eta^3q^{-1/8}\sum H_n(g)q^n$ are integral. Since $\eta^3q^{-1/8}$ is invertible in the ring
$\bbZ[[q]]$, we find that $H_n(g)\in\bbZ$ for all $n\ge 0$ and all $g\in M_{24}$.

This invertibility of $\eta^3q^{-1/8}$ directly gives us the useful implication:
\begin{equation}\label{3.2b}
lf_g\in m\bbZ[[q]]\,\Rightarrow\, lH_n(g)\equiv l\frac{w_g}{24}H_n(1)\ \ (\mathrm{mod}\ m)
\end{equation}
for all $n\ge 0$, for any choice of $l,m\in\bbZ$ and all $g\in M_{24}$.

We need to verify, for each pair $(p,K_g)\in\cP_{M24}$, that $\sum_{i=1}^k\ell_i H_n(K_i)\in p^\pi R_p$ 
for all $(\ell_1,\ldots,\ell_k)\in \cM_{(p,K_g)}$. 

A typical example is the pair $(3,2A)$. From 
Table 3 we see that $S=K_{2A}\cup K_{6A}$. Note that $\chi(2A)\equiv_3 \chi(6A)$,
for all irreducible $\chi$, so  we need to show $H_n(2A)\equiv_3
H_n(6A)$ for all $n$. 
We claim that $f_{2A}-4f_{6A}\in3\bbZ[[q]]$. To see this, first note that $f_{2A}-4f_{6A}$ is a modular form
for $\Gamma_0(6)$ with integer $q$-expansion and trivial multiplier, so by Lemma 3(a) it suffices to check that
3 divides the $q^l$-coefficient of $f_{2A}-4f_{6A}$ for all $0\le l\le 2$, which is trivial
to do from the expansions collected in Section 2. By \eqref{3.2b} or otherwise, this implies
$H_n(2A)\equiv_3 H_n(6A)$ for all $n$ and we're done.
The pairs $(2,5A),(2,7AB),(3,2B),(3,4AC),(5,1A),(5,2B),(11,1A)$ are all handled similarly.

The pair $(3,5A)$ is also easy. Here $S=5A\cup 15A\cup 15B$, and we find $\cM_{(3,5A)}
=R_3$-span$\{(0,1,-1),(1,-1,0)\}$ using obvious notation. $H_n(15A)=H_n(15B)$ for all $n$ follows from
$f_{15A}=f_{15B}$; $H_n(5A)\equiv_5 H_n(15A)$ follows from $f_{5A}\equiv_3
4f_{15A}$, proved in the usual manner using $\Gamma_0(15)$. The pairs
$(5,3A),(7,1A),(7,2A),(7,3B),(23,1A)$ are done similarly.

For $(3,7AB)$, we need $H_n(7AB)\equiv_3 H_n(21AB)$ for all $n$. Using \eqref{eq:fg}, as well as the congruences
\eqref{3.2},\eqref{E23}, this is equivalent to verifying that  $4f_{7AB}- f_{21AB}\equiv_3 2E_2^{(3)}$, 
which by Lemma 3 for $\Gamma=\Gamma_0(63)$ requires checking up to $q^{16}$. 
For $(2,3A)$, $S=3A\cup 6A\cup 12A$ and $\cM_{(2,3A)}=R_2$-span$\{(2,-2,0),(1,1,1)\}$,
so we need to verify that $H_n(3A)\equiv_2 H_n(6A)$ and $H_n(3A)+H_n(6A)
\equiv_8 2H_n(12A)$, for all $n$. Using now \eqref{E22}, this is equivalent to verifying that  $f_{3A}\equiv_2 f_{6A}$
and $f_{3A}+9f_{6A}-2f_{12A}\equiv_8 -2E_2^{(2)}$. The former requires checking up
to $q^2$ (use $\Gamma_0(6)$), while the latter requires checking up to
$q^8$ (use $\Gamma_0(24)$). The pair $(2,3B)$ is handled similarly.

The pair $(3,1A)$ has $S=1A\cup 3A\cup 3B$ and $\cM_{(3,1A)}=R_3$-span$\{(1,1,1),(0,9,3),
(0,0,9)\}$, so we need to establish $H_n(1A)\equiv_3 H_n(3B)$ and
$H_n(3A)\equiv_{27} 3H_n(3B)-2H_n(1A)$ for all $n$. These are equivalent to 
$f_{3B}\equiv_3 2E_2^{(3)}$ and $4f_{3A}-12f_{3B}\equiv_{27} 9E_2^{(3)}$.

Finally, $\cM_{(2,1A)}$ is the span of $(1,1,1,1,1,1,1)$, $(-22,-6,2,2,-2,2,0)$, 
$(44,12,28,4,4,0,0)$, $(-24,8,24,8,0,0,0)$, $(-208,-16,16,0,0,0,0)$, $(448,64,0,0,0,0,0)$, using
obvious notation. Therefore we need to verify 
\begin{align*}&H_n(4C)\equiv H_n(8A)\ \ (\mathrm{mod}\ 2)\,,\\
&H_n(4B)\equiv 2H_n(8A)-H_n(4C)\ \    (\mathrm{mod}\ 4)\,,\\
&H_n(4A)\equiv H_n(4B)+2 H_n(4C) -2H_n(8A)\ \   (\mathrm{mod}\ 8)\,,\\
&H_n(2B)\equiv 3H_n(4A)+4H_n(4B)+2H_n(4C)+8H_n(8A)\ \  (\mathrm{mod}\ 16)\,,\\ 
&H_n(2A)\equiv -H_n(2B)+4H_n(4A)+6H_n(4B)-8H_n(8A)\ \  (\mathrm{mod}\ 64)\,,\\
&H_n(1A)\equiv 7H_n(2A)-6H_n(2B)+8H_n(4A)+24H_n(4B)+32H_n(4C)-64H_n(8A)\ \ (\mathrm{mod}\ 1024)\end{align*}
for all $n$. These are equivalent to
\begin{align*}
&2f_{4C}\equiv 2f_{8A}+E_2^{(2)}\ \  (\mathrm{mod}\ 4)\,,\\  
&f_{4B}\equiv 2f_{8A}-f_{4C}\ \  (\mathrm{mod}\ 4)\,,\\  
&f_{4A}\equiv f_{4B}+2f_{4C}-2f_{8A}\ \  (\mathrm{mod}\ 8)\,,\\  
&f_{2B}\equiv 3f_{4A}+4f_{4B}+2f_{4C}+8f_{8A}+8E_2^{(2)}\ \  (\mathrm{mod}\ 16)\,,\\  
&f_{2A}\equiv -f_{2B}+4f_{4A}+6f_{4B}-8f_{8A}\ \  (\mathrm{mod}\ 64)\,,\\  
&7f_{2A}-6f_{2B}+8f_{4A}+24f_{4B}+32f_{4C}-64f_{8A}\equiv 0\ \  (\mathrm{mod}\ 1024)\,,\end{align*}  
respectively, where all of these are modular forms for $\Gamma_0(16)$.\qquad\textit{QED}\medskip

Refining this argument slightly gives us the evenness property stated in Theorem B. (Theorem
B will follow from Theorem 2 and Theorem 3).

\medskip\noindent\textbf{Theorem 2.} \textit{Each head character $H_n$ is a linear combination over
$\bbZ$ of }
\begin{eqnarray}&&\{ 2,2\rho_1,\rho_2+\overline{\rho_2}, \rho_3+\overline{\rho_3},2\rho_4,2\rho_5,2\rho_6,\rho_7+\overline{\rho_7},\rho_8+
\overline{\rho_8},2\rho_9,\nonumber\\
&&\qquad\quad\,\,\rho_{10}+\overline{\rho_{10}},2\rho_{11},2\rho_{12},2\rho_{13},2\rho_{14},2\rho_{15},2\rho_{16},2\rho_{17},2\rho_{18},2\rho_{19},2\rho_{20}\}\nonumber\end{eqnarray}

\noindent\textit{Proof.} Let $\chi$ be any $M_{24}$-class function, and $\rho\in\widehat{M_{24}}$.
Write mult$_\rho(\chi)$ for the multiplicity. Consider first $\rho$ a \textit{complex} character (i.e.  
$\overline{\rho}\not\cong\rho$). 
We know that each $f_g(\tau)\in \bbQ+q\bbZ[[q]]$. The integrality of $H_n(g)$ for each $n$
means that  mult$_\rho(H_n)=\mathrm{mult}_{\overline{\rho}}(H_n)$.

Much more difficult are the real $M_{24}$-irreps $\rho$ (i.e. $\rho\cong\overline{\rho}$). Assume
$\rho\ne 1,\rho_1$ for now. Define
\begin{equation}
M_\rho(\tau):=\sum_{n=0}^\infty \mathrm{mult}_\rho(H_n)q^n=-\sum_{K_g}|C_{M_{24}}(g)|^{-1}
\overline{\rho(g)}f_g(\tau)\,q^{1/8}\eta(\tau)^{-3}\,,\label{Mi}\end{equation}
using \eqref{eq:fg}, where the second sum is over all conjugacy classes in $M_{24}$. We know from Theorem 1
that each $M_\rho(\tau)\in q\bbZ[[q]]$ (at least for $\rho\ne 1,\rho_1$). We want to show that in fact 
$M_\rho(\tau)\in 2q\bbZ[[q]]$. Of course we can ignore the constant terms of the $f_g$ in the
following.

First note from the $M_{24}$ character table that $\rho(23A)$ and $\rho(23B)$ will be
equal and integral. Also, the centraliser $C_{M24}(23AB)$ has odd order (namely, 23).
Since $f_{23AB}(\tau)\in\bbQ+q\bbZ[[q]]$, this means we can ignore the contribution of the
classes $23AB$ to \eqref{Mi} --- together they will contribute something in $\frac{2}{23}q\bbZ[[q]]$
and thus will not affect the value of $23M_i(\tau)$ modulo 2.

Identical arguments apply to the classes $21AB$ and $15AB$. The class $11A$ can likewise
be ignored: its centraliser has odd order (namely 11), $\rho(11A)\in\bbZ$, and $f_{11A}(\tau)\in
\bbQ+2q\bbZ[[q]]$ (to see this, apply Lemma 3(a) to $\frac{1}{2}(f_{11A}+\frac{1}{6}E_2^{(11)})$). 
Finally, the classes $7AB$ and $14AB$ can all be dropped from \eqref{Mi}: they all have
the same order (namely 14) of centraliser, $\rho(7A)+\rho(7B)+\rho(14A)+\rho(14B)\in 4\bbZ$
(since $\overline{\rho}\cong\rho$), and
the $q^n$-coefficients ($n>0$) of $f_{7AB}$ and $f_{14AB}$ are integers congruent mod 2 (to
see this, apply  Lemma 3(a) to $3f_{14AB}-f_{7AB}$). Let $\widetilde{M_\rho}(\tau)$ be $23\cdot 21\cdot 5\cdot 11$ times the
sum in \eqref{Mi} restricted to the remaining classes  $K_g\in\{1A,2AB,3AB,4ABC,5A,6AB,8A,
10A,12AB\}$. Then $\widetilde{M_\rho}(\tau)\in q\bbZ[[q]]$ and $M_\rho\equiv_2 \widetilde{M_\rho}$ 
 and we need to show that $\widetilde{M_\rho}(\tau)\in 2q\bbZ[[q]]$, or equivalently
$\eta^3q^{-1/8}\widetilde{M_\rho}\in 2q\bbZ[[q]]$.

Note that $\eta^3q^{-1/8}\widetilde{M_\rho}$ is a linear combination (over $\bbQ$) of $f_g$ with
$\rho(g)\ne 0$ and $g\in\{1A,2AB,3AB,4ABC,5A,6AB,8A,10A,12AB\}$.  Let $N_\rho$ be
the least common multiple of $h_gN_g$ over those $g$,  where $h_g$ is in Table 2 and $N_g$ is 
the order of $g$.  We know from Lemma 1
 that $f_g$ is a weight-2 modular form for $\Gamma_0(h_gN_g)$ with trivial multiplier,
so $\eta^3q^{-1/8}\widetilde{M_\rho}$ is a weight-2 modular form for $\Gamma_0(N_\rho)$ with trivial
multiplier.  We collect $N_\rho$ in Table 4, together with the quantity $m_\rho/6$ where
$m_\rho=\|\mathrm{ SL}_2(\bbZ)/\Gamma_0(N_\rho)\|$. By Lemma 3, we need to show that the
first $m_\rho/6$ coefficients of $\eta^3q^{-1/8}\widetilde{M_\rho}$ are even. As in \eqref{eq:fg}
this is equivalent to showing that the first $m_\rho/6$ values of mult$_\rho(H_n)$ are even.
The evenness of these multiplicities mult$_\rho(H_n)$ has been verified for all $n\le 500$
and all real $\rho$, by Gaberdiel-Hohenegger-Volpato (private communication; see also
\cite{GHV2}), which is far more than is necessary.

The proof for $\rho=1$ and $\rho=\rho_1$ is similar, but \eqref{Mi} has to be modified.
For those $\rho$ define
\begin{align}
F_\rho(\tau)&\,:=\sum_{K_g}|C_{M_{24}}(g)|^{-1}
\overline{\rho(g)}f_g(\tau)+\frac{1}{12}E_2^{(2)}(\tau)\nonumber \\ &=-\,q^{-1/8}\eta(\tau)^{3}
\sum_{n=0}^\infty \left(\mathrm{mult}_\rho(H_n)-\frac{1}{24}H_n(1)\right)q^n+\frac{1}{12}E_2^{(2)}(\tau)\nonumber\\
&\,\equiv_2-\,q^{-1/8}\eta(\tau)^{3}
\sum_{n=0}^\infty \mathrm{mult}_\rho(H_n)\,q^n\,,\label{Frho}\end{align}
using \eqref{3.2}, \eqref{E22}, \eqref{eq:fg} and the fact that $w_g=1+\rho_1(g)$. Thus if we
can show $F_\rho(\tau)$ has even coefficients, we will be done. The expression \eqref{Frho}
together with Theorem 1 tells us $F_\rho$ has integral Fourier coefficients. Define 
$\widetilde{F}_\rho$ as above by restricting the sum to classes $K_g\in\{1A,2AB,3AB,4ABC,5A,6AB,8A,10A,12AB\}$ and multiplying by  $23\cdot 21\cdot 5\cdot 11$; as above, it suffices
to show the first 288 coefficients are even, equivalently that the multiplicities of both 1 and $\rho_1$
are even in all head characters $H_n$ for $n\le 288$. This has been done in the aforementioned
computer checks. 
 \textbf{QED to Theorem 2} 

 \medskip
\centerline {{\bf Table 4.} Data for Theorem B proof}
$${\scriptsize\vbox{\tabskip=0pt\offinterlineskip
  \def\tablerule{\noalign{\hrule}}
  \halign to 5.87in{
    \strut#&\vrule# &    
    \hfil#&\vrule#&\vrule#&\hfil#&\vrule#&    
\hfil#&\vrule#&\hfil#&\vrule#&\hfil#&\vrule#&\hfil#&\vrule#&\hfil#&\vrule#&\hfil#&\vrule#&\hfil#&\vrule#&\hfil#&\vrule#&\hfil#&\vrule#&\hfil#&\vrule#&\hfil#&\vrule#&\hfil#&\vrule#&\hfil#&\vrule#&\hfil#&\vrule#&\hfil#&\vrule#\tabskip=0pt\cr\tablerule     
&&$\,\,\rho\,\,$&\,&&\,\,1\,\,&&\,\,$\rho_1$\,&&\,\,$\rho_4$\,\,&&\,\,$\rho_5$\,\,&&\,\,$\rho_6$\,\,&&\,\,$\rho_9$\,\,&&\,\,$\rho_{11}$\,\,&&\,\,$\rho_{12}$\,\,&&\,\,$\rho_{13}$\,\,&&\,\,$\rho_{14}$\,\,&&\,\,$\rho_{15}$\,\,&&\,\,$\rho_{16}$\,\,&&\,\,$\rho_{17}$\,\,&&\,\,$\rho_{18}$\,\,&&\,\,$\rho_{19}$\,\,&&\,\,$\rho_{20}$\,\,&\cr\tablerule\tablerule 
&&$\,\,N_\rho\,$&\,&&\,\,$2^43^25$\,\,&&\,\,$2^43^25$\,\,&&\,\,$2^33\,5$\,\,&&\,\,$2^43^25$\,\,&&\,\,$2^43\,5$\,\,&&\,\,$2^43^2$\,\,&&\,\,$2^43^2$\,\,&&\,\,$2^43^25$\,\,&&\,\,$2^33^25$\,\,&&\,\,$2^43^25$\,\,&&\,\,$2^23^25$\,\,&&\,\,$2\,3^2$\,\,&&\,\,$2^43\,5$\,\,&&\,\,$2^33\,5$\,\,&&\,\,$2^33\,5$\,\,&&\,\,$2^4$\,\,&\cr\tablerule 
&&$\,\,m_\rho/6\,\,$&\,&&\,\,288\,\,&&\,\,288\,\,&&\,\,48\,\,&&\,\,288\,\,&&\,\,96\,\,&&\,\,48\,\,&&\,\,48\,\,&&\,\,288\,\,&&\,\,144\,\,&&\,\,288\,\,&&\,\,72\,\,&&\,\,6\,\,&&\,\,96\,\,&&\,\,48\,\,&&\,\,48\,\,&&\,\,48\,\,&\cr \tablerule 
 \noalign{\smallskip}}}}$$

\section{Weak Mathieu Moonshine II: Positivity}

In this section we prove that for each $n>0$ and each irreducible $\rho
\in\widehat{M_{24}}$, the multiplicities
mult$_{H_n}(\rho)$ are nonnegative.
The difficult part of this positivity proof is effectively bounding a certain series (what we call $Z_{n;h}(3/4)$ 
below) which does not
converge absolutely. There are (at least) two approaches for this: interpreting this as a
Selberg-Kloosterman zeta function as in \cite{Sel},\cite{GS}, etc; or interpreting this as
a sum over equivalence classes of quadratic forms as in Hooley \cite{Hoo}. The former method
was used by \cite{CD} to prove convergence for the Rademacher sum expressions for
the mock modular forms $q^{-1/8}\sum_{n=0}^\infty H_n(g)q^n$; the latter method (or rather
its recent reincarnation \cite{BO}) was
suggested in \cite{EHa} as a way to prove convergence for $g=1$. We need much more
than convergence: we need an explicit bound, and for this we have found the second
method more useful.

 Let's begin with an elementary observation. By the triangle inequality, 
\begin{equation}
\mathrm{mult}_{H_k}(\rho)=\sum_g\frac{H_k(g)}{|C_{M24}(g)|}\overline{\rho(g)}\ge
 \frac{H_k(1)}{|{M24}|}{\rho(1)}-\sum_{g\ne 1}\frac{|H_k(g)|}{|C_{M24}(g)|}|{\rho(g)}|
\label{posbound}\end{equation}
where the sum is over all conjugacy classes of $M_{24}$. Since there always is the trivial bound 
$\rho(1)\ge|\rho(g)|$ (which itself follows immediately from the triangle inequality), our strategy to show $\mathrm{mult}_{H_k}(\rho)>0$
is to show $H_k(1)$ is much larger in modulus than the other character values $H_k(g)$,
at least for $k$ sufficiently large.

Choose any $g\in M_{24}$. Then from Lemma 1 we read that the multiplier of 
$f_g$ is $\rho_{|g|;h_g}$, where $\rho_{n;h}$ sends $\left({a\atop c}{b\atop d}\right)\in\Gamma_0(n)$ to $e\left(\frac{ cd}{nh}\right)$. 
The multiplier $\eps$ of $\eta$ sends $\left({a\atop c}{b\atop d}\right)\in\mathrm{SL}_2(\bbZ)$ (with $c>0$)
 to
\begin{equation}
\frac{1-\I}{\sqrt{2}} \omega_{-d,c}e\left(\frac{a+d}{24c}\right)\,,\end{equation}
where
\begin{align} 
\omega_{d,c}=&\prod_{\mu=1}^k\exp(\pi\I\,\left(\left(\frac{h\mu}{k}\right)\right)\,\left(\left(\frac{\mu}{k}\right)\right)\,)
\nonumber\\ =&\left\{\begin{matrix}\left(\frac{-d}{c}\right)\,e\left(\frac{-1}{8}(c-1)+\frac{-1}{24}(c-\frac{1}{c})(
2d+d'-d^2d')\right)&\mathrm{if}\ c\ \mathrm{is\ odd}\\
\left(\frac{-c}{d}\right)\,e\left(\frac{-1}{8}(2-cd-d)+\frac{-1}{24}(c-\frac{1}{c})(
2d+d'-d^2d')\right)&\mathrm{if}\ c\ \mathrm{is\ even}\end{matrix}\right.\,.\label{sdc}\end{align}
Here, $((x))=x-\lfloor x\rfloor-1/2$ unless $x\in\bbZ$ in which case $((x))=0$.
We write here $d'$ for any solution to $dd'\equiv_c 1$.
$\left(\frac{d}{c}\right)$ denotes the Jacobi symbol. 
 The first expression for $\omega_{d,c}$ is the most familiar; the second  (due to Rademacher \cite{Rad}) is more useful for us.

We learn in eq.(6.1) of \cite{CD} (see also \cite{EHa} for the special case $g=1$) that
\begin{equation}H_k(g)=\frac{4\pi}{(8k-1)^{1/4}}\sum_{c=1}^\infty\frac{1}{|g|\,c}I_{1/2}\left(\frac{\pi}{2c\,|g|}\sqrt{8k-1}\right)\,S(k,|g|\,c,\epsilon^{-3}\rho_{|g|;h_g})\,,\label{Hk}\end{equation}
where $I_{1/2}(x)=\sqrt{\frac{2}{\pi x}}\sinh(x)$ and $S$ is (up to a constant) a generalised Kloosterman sum 
for $\Gamma_0(|g|)$:
\begin{equation}  \label{kloos}
S(k,nc,\epsilon^{-3}\rho_{n;h})=\sum_{{0<d\le nc\atop \mathrm{gcd}(d,nc)=1}}
\omega_{d,nc}^{-3}e\left(\frac{-cd}{h}\right)\,e\left(\frac{kd }{nc}\right)\end{equation}


Fix integers $n,h,k>0$ and define $L=\frac{\pi}{n}\sqrt{2k-1}$ (so $\pi \sqrt{2k-1}/(2nc)<1$
if $c>L$).  We have the elementary bounds: 
\begin{align}
&\left|I_{1/2}(x)-\sqrt{\frac{2x}{\pi}}\right|\le\frac{1}{5}\sqrt{\frac{2x^5}{\pi}}\qquad\mathrm{for}\
0<x<1\,;\label{smallx}\\
&|I_{1/2}(x)|<\frac{e^x}{\sqrt{2\pi x}}\qquad\mathrm{for\ all}\ 0<x\,;\label{bigx}\\
&|S(k,nc,\epsilon^{-3}\rho_{n;h})|\le 1\qquad\mathrm{for\ all}\ c\in\bbZ_{>0}\,\label{naive}\end{align}
(in fact Lemma 6 below will give us the Weil bound $|S|\le O(\sqrt{c})$, but \eqref{naive}
is adequate). These inequalities immediately imply  the crude bounds
\begin{align}
&\left|\frac{1}{(8k-1)^{1/4}}\sum^{\lfloor L\rfloor}_{c=2}\frac{1}{c}I_{1/2}\left(\frac{\pi\sqrt{8k-1}}{2nc}\right)
\,S(k,nc,\epsilon^{-3}\rho_{n;h})\right|\nonumber\\
&\ \ \ \ \le \frac{\sqrt{n}}{\pi\sqrt{8k-1}}\sum_{c=2}^{\lfloor L\rfloor}
\frac{1}{\sqrt{c}}e^{\pi\sqrt{8k-1}/(2nc)}<\frac{1}{\sqrt{2n}}e^{\pi\sqrt{8k-1}/(4n)}\,,\\
&\left|\frac{1}{(8k-1)^{1/4}}\sum_{c=\lceil L\rceil}^\infty \frac{1}{nc}I_{1/2}\left(\frac{\pi\sqrt{8k-1}}{2nc}\right)
\,S(k,nc,\epsilon^{-3}\rho_{n;h})\right|\nonumber\\
&\ \ \ \ \le \left|\sum_{c=\lceil L\rceil}^\infty (nc)^{-3/2}S(k,nc,\epsilon^{-3}\rho_{n;h})\right|+\frac{\pi^2(8k-1)}{20n^{7/2}}\sum_{c=\lceil L\rceil}^\infty c^{-7/2}\nonumber\\
&\ \ \ \
\le |Z_{n;h}(3/4)|+\frac{\pi\sqrt{8k-1}}{n^{5/2}}+\frac{7\pi^2(8k-1)}{80n^{7/2}}\,,\end{align}
where $Z_{n;h}$ is the Selberg-Kloosterman zeta function 
\begin{equation} \label{zeta}
Z_{n;h}(s)=\sum_{c=1}^\infty S(k,nc,\epsilon^{-3}\rho_{n;h})\,(nc)^{-2s}\,.\end{equation}

Thus it suffices to find an effective bound for a Selberg-Kloosterman zeta function
at $s=3/4$.
We know from \cite{CD} that the defining series \eqref{zeta} converges there, but it follows from
Lemma 4 below that this convergence is not absolute, and in any event it seems difficult to use
the analysis of \cite{CD} to obtain an explicit bound  (as a function in $k$). 
The key first step in identifying such an effective bound is to rewrite these generalised Kloosterman sums more sparsely.  Our calculation resembles that of
\cite{Wh}; a more elegant approach attributed to Selberg and worked out by Rademacher
\cite{Rad1} is available but we were unable to generalise it to our context.
Incidentally, Lemma 4 should imply that these $S$'s have some multiplicative
properties (hence their $Z$'s have Euler-like products).

\medskip\noindent\textbf{Lemma 4.} \textit{Let $c\in\bbZ_{>0}$, $k\in\bbZ$, $h|n$. Then}
\begin{equation}
S(k,nc,\epsilon^{-3}\rho_{n;h})=\frac{-\I
\sqrt{nc}}{2}\sum_{0\le m<4nc\atop m^2\equiv_{8nc}1-8k+8c^2n/h}e\left(\frac{m}{4nc}\right)\,.
\end{equation}

\noindent\textit{Proof.} We will begin with the proof of $n=h=1$. Writing $m=2\ell+1$, define
\begin{equation}
B_c(k):=\frac{-\I
\sqrt{c}}{2}\sum_{m=0\atop m^2\equiv_{8c}1-8k}^{4c-1}(-1)^{(m-1)/2}e\left(\frac{m}{4c}\right)
=\frac{-\I\sqrt{c}}{2}e\left(\frac{1}{4c}\right)\sum^{2c-1}_{\ell=0\atop \ell^2+\ell\equiv_{2c}-2k}(-1)^\ell
e\left(\frac{\ell}{2c}\right)\,.\nonumber\end{equation}
As we manifestly have $B_c(k+c)=B_c(k)$ (i.e. $B_c$ is a class function for $\bbZ_c$), we can formally write it as a combination of irreducible $\bbZ_c$-characters:
\begin{equation}
B_c(k)=\sum_{d=0}^{c-1}R_{d,c}\,e\left(\frac{dk}{c}\right)\,,\end{equation}
for coefficients 
\begin{align}\label{Rdc}
R_{d,c}=\frac{-\I}{2\sqrt{c}}e\left(\frac{1}{4c}\right)\sum_{j=0}^{c-1}e\left(\frac{-dj}{c}\right)
\sum^{2c-1}_{\ell=0\atop \ell^2+\ell\equiv_{2c}-2k}(-1)^\ell
e\left(\frac{\ell}{2c}\right)&\nonumber\\
=\frac{-\I}{2\sqrt{c}}e\left(\frac{1}{4c}\right)\sum_{\ell=0}^{2c-1}(-1)^\ell
e\left(\frac{d\,(\ell^2+\ell)}{2c}+\frac{\ell}{2c}\right)\,.&\end{align}
It is elementary to show that the sum on the right-side of \eqref{Rdc} vanishes when
gcd$(c,d)>1$: write $m=\mathrm{gcd}(c,d)$ and $\ell=s2c/m+r$, so $\sum_{\ell=0}^{2c-1}=
\sum_{r=0}^{2c/m-1}\sum_{s=0}^{m-1}$, and notice that $\sum_s=0$ for each $r$ (when $m>1$).

Thus we can restrict to gcd$(c,d)=1$. We want to show $R_{d,c}=\omega_{d,c}^{-3}$. 
Choose $d'\in\bbZ$ so that $1=e\left(\frac{(dd'-1)(c+1)}{2c}\right)$; this permits us to rewrite
\eqref{Rdc} as
\begin{equation}
R_{d,c}=\frac{-\I e({1}/{4c})}{2\sqrt{c}}G(d,d\gamma;2c)\end{equation}
where $\gamma:=d'c+d'+1$ and $G(a,b;c)$ is the generalised Gauss sum $G(a,b;c)=\sum_{\ell=0}^ce((a\ell^2+b\ell)/c)$.
Note that $\gamma:=d'c+d'+1$ is even iff $c$ is even. 
When $c$ is even, we can complete squares and obtain
\begin{align}\label{even}
R_{d,c}&=\frac{-\I \,e({1}/{4c})}{2\sqrt{c}}\,e\left(\frac{-d\gamma^2}{8c}\right)\,G(d,0;2c)\nonumber\\
&=
\left(\frac{2c}{d}\right)\,e\left(\frac{-(d-1)^2}{16}+\frac{2-c-d\gamma^2}{8c}\right)\,.\end{align}
When $c$ is odd, use $\frac{1}{2c}=\frac{-1}{2}+\frac{(c+1)/2}{c}$ to write $G(d,d\gamma;2c)
=G(-d,-d\gamma;2)\,G(d(c+1)/2,d\gamma(c+1)/2;c)$; the left generalised Gauss sum equals
2, while the right is evaluated by completing squares as usual, and we obtain (for $c$ odd)
\begin{equation}\label{odd}
R_{d,c}=\left(
\frac{d(c+1)/2}{c}\right)\,e\left(\frac{(c-1)^2}{16}+\frac{2-2c-d(c+1)^3(d'+1)^2}{8c}\right)\,.\end{equation}

Now, for $m$ odd we have $\left(\frac{2}{m}\right)=(-1)^{(m^2-1)/8}$ and $\left(\frac{-1}{m}\right)
=(-1)^{(m-1)/2}$. Consider first $c$ even. Then directly from \eqref{sdc} and \eqref{even}
we obtain
\begin{equation}\frac{\omega^{-3}_{d,c}}{R_{d,c}}=e\left(\frac{-2-2c-d-d'-2cd+2dd'-c^2d+c^2d'+
d^2d'+dd'{}^2+2cdd'{}^2+c^2dd'{}^2-c^2d^2d'}{8c}\right)\,.\label{ratio-ev}\end{equation}
Because $dd'\equiv_{2c}1$ when $c$ is even we can define an integer $L$ by $dd'=1+2cL$,
and we know that $d\equiv_4d'$. Then \eqref{ratio-ev} collapses to
\begin{equation}\frac{\omega^{-3}_{d,c}}{R_{d,c}}=e\left(\frac{dL+d'L+2L}{4}\right)=1\,,\end{equation}
as desired. The proof for $c$ odd is similar: the Jacobi symbol $\left(
\frac{d(c+1)/2}{c}\right)$ equals $\left(\frac{2d}{c}\right)$; define $L$ by $dd'=1+cL$
and use the congruences $c^3\equiv_{8c} c$ and $4c^2\equiv_{8c}4c$.

Replacing $c$ everywhere with $nc$ establishes the $h=1$ case of Lemma 4. Arbitrary $h|n$ is handled through the elementary 
observation that $S(k,nc,\epsilon^{-3}\rho_{n;h})=S(k-c^2n/h,nc,\epsilon^{-3}\rho_{n;1})$.
 \qquad \textbf{QED to Lemma 4}\medskip

Hooley's method is based on the $n=1$ case of Lemma 5 below; it was more recently modified
slightly by \cite{BO}. We need a more serious revision. A binary quadratic form $Q(x,y)
=\alpha x^2+\beta xy+\gamma y^2=:[\alpha,\beta,\gamma]$ is called \textit{integral} if $\alpha,\beta,\gamma\in\bbZ$, and \textit{positive-definite} if $\alpha$ and $\gamma$ are both positive. (See e.g. Chapter 12 of \cite{Hua} for a rather complete introduction to the basic theory
of quadratic forms, as is relevant here.) The discriminant is $\beta^2
-4\alpha\gamma$.  For any $C\in\bbZ_{>0}$ and $D\in\bbZ_{<0}$, let $\mathcal{Q}(C,D)$ denote the set of all
triples $(Q;r,s)$, where $Q$ is integral and positive-definite with discriminant $D$ and
where $r,s$ are coprime integers satisfying $Q(r,s)=C$. Any $\left(\begin{matrix} a&b\\ c&d\end{matrix}\right)\in\mathrm{SL}_2(\bbZ)$ acts on $(Q;r,s)\in\mathcal{Q}(C,D)$ by sending
$Q(x,y)$ to  $Q(ax+by,cx+dy)$ and  $(r,s)$ to $(dr-bs,-cr+as)$. In particular, it is easy to
verify that if $(Q;r,s)$ and $(Q';\tilde{r},\tilde{s})$ 
are in the same SL$_2(\bbZ)$-orbit then
$Q,Q'$ have the same discriminant and $Q(r,s)=Q'(\tilde{r},\tilde{s})$. Then Hooley
observes that there is a bijection between the integers $0\le m<2C$ satisfying $m^2\equiv_{4C} 1$,
and SL$_2(\bbZ)$-orbits in $\mathcal{Q}(C,D)$.

We need to generalise this in two ways. First, choose any integer $n\ge 1$. Write
$\mathcal{Q}_n(C,D)$ denote the set of all triples $(Q;r,ns)$, where $Q=[n\alpha,\beta,\gamma]$
is positive-definite and of discriminant $D$, $\alpha,\beta,\gamma,r,s\in\bbZ$, gcd$(r,ns)=1$,
and $Q(r,ns)=nC$. For example, $\mathcal{Q}(C,D)=\mathcal{Q}_n(C,D)$.
It is elementary to show that the group $\Gamma_0(n)$ acts on $\mathcal{Q}_n(C,D)$. 

Secondly, suppose  that $h$ divides gcd$(n,24)$ and that gcd($n/h-1,h)=1$.  Write $n'=n/h$. Let
$\mathcal{Q}_{n;h}(C,D)$ denote the set of all triples $(Q;r,ns)$ where $Q=[n\alpha,\beta,\gamma/h]$ is positive-definite and of discriminant $D$, $\alpha,\beta,\gamma,r,s\in \bbZ$,
$\gamma\equiv_h\alpha$, gcd$(r,n's)=1$, and $Q(r,ns)=C$. Of course $\mathcal{Q}_{n;1}(C,D)
=\mathcal{Q}_n(C,D)$. Recall the group $\Gamma_0(n|h)$, which can be defined equivalently
as either the conjugate of $\Gamma_0(n')$ by $\left({h\atop 0}{0\atop 1}\right)$,
or as the set of all determinant-1 matrices of the form $\left({a\atop nc}{b/h\atop d}\right)$
for $a,b,c,d\in\bbZ$. Let  $\Gamma_0(n;h)$ denote the set of all  $\left({a\atop nc}{b/h\atop d}\right)
\in\Gamma_0(n|h)$ for which $ac\equiv_h bd$. Then it is easy to verify $\Gamma_0(n;h)$
is a group, and we see shortly that it
acts on $\mathcal{Q}_{n;h}(C,D)$.

The reason 24
arises here (and elsewhere) is because, for any divisor $d$ of 24, any integer $\ell$ coprime 
to $d$ satisfies $\ell^2\equiv_d1$.
The defining condition $ac\equiv_h bd$ of $\Gamma_0(n;h)$ is equivalent to requiring that there 
is  some $\ell$ (depending on $a,b,c,d$) coprime to $h$ 
with $a\equiv_h\ell d$ and $b\equiv_h \ell$. To see this equivalence, run through each prime power $p^\nu$
exactly dividing $h$; the determinant condition $ad-n'bc=1$ tells us that either $a$ and
$d$ are both coprime to $h$ (in which case require $\ell\equiv ad$ (mod $p^\nu$)),
or $b,c$ are both coprime to $h$ (in which case take $\ell\equiv bc$ (mod $p^\nu$)). 
Then $a\equiv_{p^\nu} \ell d$ resp.\ $b\equiv_{p^\nu} \ell c$, and the other congruence comes
from $ac\equiv_{p^\nu} bd$. Running through all $p$, we obtain an $\ell$ defined mod $h$
which has all the desired properties.

\medskip\noindent\textbf{Lemma 5(a)} \textit{Let $D\in\bbZ_{<0}$, $C\in\bbZ_{>0}$. There is
a one-to-one correspondence between the set of integers $m$, $0\le m<2nC$, satisfying
$m^2\equiv D$ (mod $4nC$), and $\Gamma_0(n)$-orbits in $\mathcal{Q}_n(C,D)$.}\medskip 

\noindent\textbf{(b)} \textit{There is
a one-to-one correspondence between the set of integers $m$, $0\le m<2nC$, satisfying
$m^2\equiv D+4C^2n'$ (mod $4nC$), and $\Gamma_0(n;h)$-orbits in $\mathcal{Q}_{n;h}(C,D)$.}\medskip

\noindent\textit{Proof.} Let's begin with the simpler part (a).  First note that there is an obvious bijection between each $m\in\bbZ$
with $m^2\equiv_{4nC}D$, and each quadratic form  $Q=[nC,m,\gamma]$ (necessarily positive-definite)
of discriminant $D$,   where $\gamma=(m^2-D)/(4nC)$. 

Given any $(Q;r,ns)\in\mathcal{Q}_n(C,D)$, $Q=[n\alpha,\beta,\gamma]$, choose any $\tilde{r},\tilde{s}$
such that $\left(\begin{matrix} r&\tilde{r}\\ ns&\tilde{s}\end{matrix}\right)\in\Gamma_0(n)$ (this
is possible since gcd$(r,ns)=1$), and
define the \textit{root} $m:=2n\alpha r\tilde{r}+\beta\,(r\tilde{s}+ns\tilde{r})+2\gamma ns\tilde{s}$. Then $
\left(\begin{matrix} r&\tilde{r}\\ ns&\tilde{s}\end{matrix}\right)\in\Gamma_0(n)$ sends $(Q;r,ns)$ to 
$([nC,m,\gamma'];1,0)$
where $\gamma'=(m^2-D)/(4nC)$. The other choices of $\tilde{r},\tilde{s}$ are $\tilde{r}+Lr,\tilde{s}+Lns$ for any $L\in
\bbZ$, which correspond to roots $m+L2nC$, and so the root taken mod $2nC$ is a well-defined
function $m(Q;r,ns)\in\bbZ_{2nC}$ on $\mathcal{Q}_n(C,D)$. The desired bijection in part (a)
is this root map (mod $2nC$). 

Now turn to part (b).  Here, there is an elementary bijection between
each $m\in\bbZ$
with $m^2\equiv_{4nC}D+4n'C^2$, and each positive-definite quadratic form  $Q=[nC,m,\gamma/h]$ 
of discriminant $D$,   where $\gamma\equiv_h C$. Choose any $([n\alpha,\beta,\gamma/h];
r,ns)\in\mathcal{Q}_{n;h}(C,D)$. Then $\left({a\atop nc}{b/f\atop d}\right)\in\Gamma_0(n|h)$
will send $(r,ns)$ to $(ra+n'sb,nrc+nsd)$ (the desired form), and $[n\alpha,\beta,\gamma/h]$ to
\begin{equation}
[n(\alpha a^2+\beta ac+n'\gamma c^2),2n'\alpha ab+\beta\,(ad+n'bc)+2\gamma n'cd,
(n'\alpha b^2+\beta bd+\gamma d^2)/h]\,.\label{qf}\end{equation}
Then \eqref{qf} will lie in $\mathcal{Q}_{n;h}(C,D)$ if
\begin{equation}
a^2\equiv_h d^2\,,\ \ ac\equiv_hbd\,,\ \ n'c^2\equiv_h n'b^2\,.\label{congruences}\end{equation}
As explained above, the condition $ac\equiv_hbd$ is equivalent to the existence of an
$\ell$ coprime to $h$ satisfying $a\equiv_h\ell d$ and $b\equiv_h\ell c$, and such an
$\ell$ forces the other two congruences to be satisfied. In other words, the matrices in
$\Gamma_0(n|h)$ satisfying \eqref{congruences} form the group $\Gamma_0(n;h)$. 

Since $r$ and $n's$ are coprime, we can find integers $\tilde{r},\tilde{s}$ such that
$r\tilde{s}-n's\tilde{r}=1$. We claim there is some $L\in\bbZ$, unique modulo $h$,
 such that $\left({r\atop ns}{(Lr+\tilde{r})/h
\atop Ln's+\tilde{s}}\right)\in\Gamma_0(n;h)$. To see this, choose any prime power $p^\nu$
exactly dividing $h$. If $r$ is coprime to $p$ choose $\ell\equiv_{p^\nu}1-n' s^2$ and
$L\equiv_{p^\nu} r'\,(-\tilde{r}+\ell s)$ where $r'r\equiv_{p^\nu}1$; otherwise, $n's$ will be
 coprime to $p$, so choose $s'$ by $s'n's\equiv_{p^\nu}1$, $\ell_{p^\nu}r^2-n'$, and $L\equiv_{p^\nu}s'\,(-\tilde{s}+\ell r)$. We are using here that $1-n'$ is coprime to $h$ and hence $p$.
 Then $\ell$ and $L$ are defined mod $h$ by running through all $p$.
To see uniqueness of $L$ modulo $h$, note that $ \left({r\atop ns}{(Lr+\tilde{r})/h
\atop Ln's+\tilde{s}}\right) \left({r\atop ns}{(L'r+\tilde{r})/h
\atop L'n's+\tilde{s}}\right)^{-1}= \left({1\atop 0}{(L-L')/h
\atop 1}\right)\in\Gamma_0(n;h)$ forces $(L-L')/h\in\bbZ$ by definition of $\Gamma_0(n;h)$.

So we may assume $\left({r\atop ns}{\tilde{r}/h
\atop \tilde{s}}\right)\in\Gamma_0(n;h)$. Define the root 
$m:=2n'\alpha r\tilde{r}+\beta\,(r\tilde{s}+n's\tilde{r})+2\gamma n's\tilde{s}$ as before. Then $
\left(\begin{matrix} r&\tilde{r}/h\\ ns&\tilde{s}\end{matrix}\right)$ sends $(Q;r,ns)$ to 
$([nC,m,\gamma'];1,0)\in\mathcal{Q}_{n;d}(C,D)$
where $\gamma'=(m^2-D)/(4nC)$. From the uniqueness of the previous paragraph, we have that 
$\left({r\atop ns}{(Lr+\tilde{r})/h
\atop Ln's+\tilde{s}}\right)\in\Gamma_0(n;h)$ iff $L\in h\bbZ$, which corresponds to root
 $m+L2nC$, and so again the root taken mod $2nC$ is a well-defined
function $m(Q;r,ns)\in\bbZ_{2nC}$ on $\mathcal{Q}_{n;d}(C,D)$. \qquad \textbf{QED to Lemma 5}\medskip

We are interested in $D=1-8k$, $n=2n_g$, $h=h_g$ and $C=c$. Restrict attention here to $k\ge 1$ (not a problem since we are only interested in large $k$), so $D<0$. In fact, to sharpen
slightly some of our bounds, we'll choose $k\ge 5$. Continue to write
$n'=n/h$.

Note that if both $(Q;r,ns),(Q;\tilde{r},n\tilde{s})$ (same $Q$) lie in $\mathcal{Q}_{n;k}(C,D)$, then there is an
automorphism $g\in\Gamma_0(n;k)$ of $Q$ sending $(r,ns)$ to $(\tilde{r},n\tilde{s})$. For $D<0$
(the case we are interested in), each such automorphism $g$ must have finite order (since completing
squares in $\alpha=\alpha a^2+\beta ac+n'\gamma c^2$ and $\gamma=n'\alpha b^2+\beta bd
+\gamma d^2$ bounds the matrix entries $a,b,c,d$ of $g$). Now, $\Gamma_0(n|h)$ (which
contains $\Gamma_0(n;h)$) is conjugate to $\Gamma_0(n')$, so all of its elements of finite
order have orders 2, 4, 6. In fact for us, the stabiliser stab$([n\alpha,\beta,\gamma/h];\Gamma_0
(n;h))$  of $[n\alpha,\beta,\gamma/h]$ in $\Gamma_0(n;h)$ will always be $\pm 1$. To see this,  
first note stab$([n\alpha,\beta,\gamma/h];\Gamma_0(n;h))$ is a subgroup of stab$([n\alpha,\beta,
\gamma/h];\Gamma_0(n|h))$, which is isomorphic to stab$([n'\alpha,\beta,\gamma];\Gamma_0
(n;h))$, which is in turn a subgroup of stab$([n'\alpha,\beta,\gamma];\mathrm{SL}_2(\bbZ))$.
For a \textit{primitive} form $[n'\alpha,\beta,\gamma]$ (i.e. when $n'\alpha,\beta,
\gamma$ don't have a common factor), nontrivial stabilisers occur only for discriminant $D=-3,-4$.
But our discriminant satisfies $D\equiv 1$ (mod 8), so even if  $[n'\alpha,\beta,\gamma]$ is
imprimitive, it can never have a nontrivial stabiliser in SL$_2(\bbZ)$. 

We also need a bound on $\Gamma_0(n;h)$-equivalence class representatives $[n\alpha,\beta,
\gamma/h]$. For this purpose, observe that $\Gamma_0(n;h)$ contains $\Gamma_0(nh)$, which has finite index
$nh\prod_p(1+1/p)$ in SL$_2(\bbZ)$. The cosets
for $\Gamma_0(N)\backslash\mathrm{SL}_2(\bbZ)$ are in bijection with pairs $(c,d)\in
\bbZ_{>0}^2$ where $c|N$ and $1\le d\le N/c$ satisfies gcd$(d,c,N/c)=1$. To any such pair
$(c,d)$, a representative of that coset is $\left({a\atop c}{b\atop d'}\right)$ where $d'$ is
coprime to $c$ and satisfies $d'\equiv_{N/c}d$, and $a,b$ are any  integers satisfying
$ad'-bc=1$. Now, $c=N$ corresponds to the identity coset, so for it take $c=0$ instead; then in all
cases we have $c\le N/2$. We can 
always choose $|a|\le c/2$ by adjusting $b$ appropriately (at least when $c\ne 0$), hence $|a|\le N/4$ (true even for $c=0$
unless $N<4$).
Now, over SL$_2(\bbZ)$, any positive-definite discriminant $D$ quadratic form (with coefficients
in $h^{-1}\bbZ$) is equivalent to some $[\alpha',\beta',\gamma']$
(namely Gauss' reduced form) with  $|\beta'|\le\alpha'\le\gamma'$, where 
$\alpha'\le\sqrt{|D|/3}$ and $\gamma'\le 11h|D|/39$. Combining with \eqref{qf}, and noting that
we can always force $|\beta|\le n\alpha$ by applying multiples of $\left({1\atop 0}{1\atop 1}\right)
\in\Gamma_0(nh)$,
we see that any such quadratic form is equivalent over $\Gamma_0(nh)$, hence over
$\Gamma_0(n;h)$, to some $[n\alpha,\beta,\gamma/h]$ satisfying the (crude but adequate)
bounds 
\begin{align}
&0<\alpha\le{n^2h^2\over 16n}\left( 3\sqrt{{|D|\over 3}}+44h|D|/39\right)<\frac{1}{12}{nh^3}\,|D|\,,\\
&\ \ \  \  \ |\beta|
<n^2h^3\,|D|/12\,,\\ &0<\gamma\le
(|D|+\beta^2)/(4n'\alpha)<\left\{\begin{matrix}|D|/7&\mathrm{if}\ n=2\\ n^2h^4\,|D|/44&\mathrm{otherwise}\end{matrix}\right.\,.
\label{qfbound}\end{align}
(The given bound on $\alpha$, hence $\beta$, is also true for the identity coset $c=0$.) 
Implicit in these derivations are $2h\le n$, $n$ even, and $|D|\ge 39$.

The point is that the series $Z_{n;h}(3/4)$ can be rewritten as a sum over the finitely many 
$\Gamma_0(n;h)$-orbit representatives
$[n\alpha,\beta,\gamma/h]$, and over integers $r,s$. Over-counting by a factor
of 2, we can take this to be a free sum over $r,s$.

\medskip\noindent\textbf{Theorem 3.} \textit{Choose any $k,n\in\bbZ_{>0}$, $k\ge 5$, and
any $h$ dividing gcd$(n,24)$ such that gcd($n/h-1,2h)=1$.  Then for $n\ne 2$,}
\begin{align}
\left|Z_{n/2;h}(3/4)\right|&\le \left(\prod_{p|(nh)}\frac{p+1}{p}\right)\left(1+2.13|D|^{1/8}\log\,|D|\right)
\times\nonumber\\ &
\left((6.124n^{35/6}h^{47/6}
-3.09n^{23/4}h^{31/4}+64.32n^{29/6}h^7-23n^{19/4}h^7){|D|}\right.\nonumber\\
&\left.+(.146n^{47/6}h^{65/6}-.114n^{31/4}h^{43/4}+2.51n^{35/6}h^{10}-.74n^{23/4}h^{10})|D|^{3/2}\right)\,.\nonumber
\end{align}
\textit{where we write $D=1-8k$. For $n=2$, this bound on the right-side should be replaced with
$(3872|D|+213|D|^{3/2})(1+2.13|D|^{1/8}\log|D|)$.}

\medskip\noindent\textit{Proof.} Note the $n$ in the statement of Theorem 3 agrees with
that of Lemma 5, but is twice that of Lemma 4.
We need to bound the limit as $X\rightarrow \infty$ of
\begin{equation}
Z_{n/2;h}(s;X)=\sum_{c=1}^{X}S(k,nc,\epsilon^{-3}\rho_{n/2;h})\,(nc/2)^{-2s}\,.\end{equation}
Thanks to Lemmas 4 and 5, we can write its value at $s=3/4$  as 
\begin{equation}\label{ZnhX}
\frac{-\I}{2}\sum_{[n\alpha,\beta,\gamma/h]}\sum_{c(\alpha,\beta,\gamma;r,s)\le X
}\!\!\!\!\!\!\!\!\!{}'\,\,\,\,\,\,\frac{(-1)^{(m(\alpha,\beta,\gamma;r,s)-1)/2}}{nc(\alpha,\beta,\gamma;r,s)}\,e\left(\frac{m(\alpha,\beta,\gamma;r,s)}{2n\,c(\alpha,\beta,\gamma;r,s)}\right)\,,\end{equation}
where the first sum is over representatives of the $\Gamma_0(n;h)$-orbits of  positive-definite quadratic forms of discriminant
$D=1-8k$, and the second sum is over all  integers $r,s$ satisfying {gcd}$(r,n's)=1$ and
the given inequality (the prime on the sum denotes that coprime condition). The quantities
$m,c$ are defined by
\begin{equation}m(\alpha,\beta,\gamma;r,s)=2n'\alpha r\tilde{r}+\beta\,(r\tilde{s}+n's\tilde{r})+2\gamma n's\tilde{s}\,,\
c(\alpha,\beta,\gamma;r,s)=\alpha r^2
+\beta rs+n'\gamma s^2\,,\nonumber\end{equation}
where $\tilde{r},\tilde{s}$ are any integers satisfying $\left({r\atop ns}{\tilde{r}/h\atop \tilde{s}}\right)\in
\Gamma_0(n;h)$ --- which pair is chosen won't affect the value of that summand.
The factor of 2 in \eqref{ZnhX} comes from the aforementioned redundancy that 
different $(r,ns)$ can lie in the same orbit.

Write this inner sum $\sum_{r,s}$ (for each choice of $\alpha,\beta,\gamma$) as $\sum_>+\sum_<$, depending on whether or not
$|r|>n|s|$ ($|r|=n|s|$ would contradict gcd$(r,s)=1$) .
We will bound $Z_{n/2;h}(3/4;X)$ by considering separately the contributions of $\sum_>,\sum_<$. 

The arguments for $\sum_>$ and $\sum_<$ are completely analogous, and so we will 
consider in detail only $\sum_>$. Since then $r>0$, we can write 
\begin{equation}\label{separ}
{m\over 2nc}=\frac{2n'\tilde{r}\,(\alpha r^2+\beta rs+n'\gamma s^2)+2n'\gamma s+\beta r}{2nr\,(\alpha r^2+\beta rs
+n'\gamma s^2)}\equiv_1-\frac{us'}{R}+\frac{2n'\gamma s+\beta r}{2nr\,(\alpha r^2+\beta rs
+n'\gamma s^2)}\end{equation}
where $\delta=\mathrm{gcd}(s,h)$, $h'=h/\delta$, $R=rh'$, $s'\in\bbZ$ is any inverse of $s/\delta$ 
mod $R$, and $u\in\bbZ$ is defined mod $R$ by $u\equiv_r-\delta^{-1}(n')^{-1}$, $u\equiv_{h'}1-s^2n'$, and (if $3|\mathrm{gcd}(r,h)$) $u\equiv_9-(n')^{-1}$
(note that gcd$(r,h)\ne 1$ implies it equals 3, in which case $s$ is coprime to 3 and
$n'\equiv_3-1$). 
Writing 
\begin{equation}\alpha r^2+\beta rs+n'\gamma s^2=n'\gamma\,(s+\frac{\beta r}{2n'\gamma})^2
-\frac{Dr^2}{4n'\gamma}=\alpha\,(r+\frac{\beta s}{2\alpha})^2-\frac{Ds^2}{4\alpha}\,,\label{compsq}\end{equation}
 the first inequality gives
\begin{equation}|r|\le \sqrt{\frac{X4n'\gamma}{|D|}}=C\sqrt{X}\,,\ \ \ \ 
S_-\le s\le S_+\,,\end{equation}
where $C=\sqrt{4n'\gamma/|D|}$ and
\begin{align}
&S_+=\mathrm{min}\left\{ \frac{|r|}{n},-\frac{\beta r}{2n'\gamma}+\frac{1}{2n'\gamma}\sqrt{4n'\gamma X+{Dr^2}}\right\} \,,\\
&S_-=\mathrm{max}\left\{-\frac{|r|}{n},-\frac{\beta r}{2n'\gamma}-\frac{1}{2n'\gamma}\sqrt{4n'\gamma X+{Dr^2}}\right\}\,.\end{align}

We compute
\begin{equation}(-1)^{m}=
(-1)^{n'\,(\alpha r\tilde{r}+\beta s\tilde{r}+\gamma s\tilde{s})+(\beta-1)/2}=(-1)^{(\beta-1)/2}
\end{equation}
using the evenness of $n'$ and the determinant relation $r\tilde{s}-n's\tilde{r}=1$ ($\beta$ is
odd because of the discriminant condition $\beta^2-4n'\alpha\gamma=1-8k$). Now put
\begin{equation}
\varphi(r,s)=\frac{(-1)^{(\beta-1)/2}}{\alpha r^2
+\beta rs+\gamma n's^2}\, e\left(\frac{2n'\gamma s+\beta r}{2nr\,(\alpha r^2+\beta rs
+\gamma n's^2)}\right)\,;\end{equation}
we need to bound
\begin{equation}
\sum_>=\sum_{|r|\le C\sqrt{X}}\,\sum_{\delta|h}\sum_{S_-\le s\le S_+\atop \mathrm{gcd}(s,h)=\delta}\!\!\!\!\!\!{}'\,\,
e\left(\frac{us'}{hr/\delta}\right)\varphi(r,s)
=\sum_{|r|\le C\sqrt{X}}\,\sum_{\delta|h}\sum_{S_-/\delta\le S\le S_+/\delta\atop \mathrm{gcd}(S,R)=1}\!\!\!\!\!\!\!
e\left(\frac{uS'}{R}\right)\varphi(r,\delta S)
\,,\label{sum11}\end{equation}
where we use \eqref{separ}, we write $R=hr/\delta$ as before and $S'=s'$ is any inverse of $S$ mod $R$.
Rewrite the sum over $S$ in \eqref{sum11} (for fixed $\alpha,\beta,\gamma,r,\delta$) using partial summation (the discrete analogue
of integration by parts):
\begin{align}
&\sum_{S_-/\delta\le S\le S_+/\delta}\!\!\!\!\!\!\!\!\!{}'\,\,\,\, e\left(\frac{uS'}{R}\right)\varphi(r,\delta S)=\!\!\!\sum_{S_-/\delta\le \sigma\le S_+/\delta}\!\!\!\!\!\!\!
g(\sigma)\,(\varphi(r,\delta\sigma)-\varphi(r,\delta\sigma+\delta))\nonumber\\ &\ \ \qquad \qquad \ \ \ \ \ \ +g(\lfloor S_+/\delta\rfloor)\,\varphi(r,\delta\lfloor S_+/\delta\rfloor)\,,\label{parsum}
\end{align}
where
\begin{equation}
g(\sigma)=\sum_{S_-/\delta\le S\le \sigma\atop \mathrm{gcd}(S,R)=1}\!\!\!\!\!\!'\,e(uS'/R)\end{equation}
denotes the \textit{incomplete} Kloosterman sums ($R$ is implicit). We have (for $n>2$)
\begin{align}
&|\varphi(r,s)-\varphi(r,s+\delta)|\le \frac{16n'{}^2\gamma^2\delta}{D^2r^4}\left|\frac{|r|}{2}\sqrt{|D|}+n'\gamma\delta\right|+\frac{2\pi}{2n|r|}\frac{(4n'\gamma)^2}{|D|r^2}\frac{2}{2|r|\sqrt{|D|}}
\nonumber\\&\qquad\qquad \qquad \ \ \ \ \ \ \ \ \ \ \le .0042\frac{n^6h^6}{|r|^3}\delta\sqrt{|D|}+.000188\frac{n^9h^9}{r^4}\delta^2|D|+.0083\frac{n^{5}h^6}{r^4}\sqrt{|D|} \,,\label{deltaphi}\\
&\ \ \ \ \ \ \ \ \ |\varphi(r,\delta\lfloor S_+/\delta\rfloor)|\le {4n'\gamma}/({|D|r^2})\le .091n^3h^3r^{-2}\,,\label{phi}\end{align}
where we use repeatedly \eqref{compsq}, as
well as Taylor's remainder $|e(\theta)-1|\le 2\pi\theta_0$ for $|\theta|\le \theta_0$,
and the elementary bound $x/(x^2+a^2)\le 1/(2|a|2)$. Thus
all that remains is to bound $g(\sigma)$. The order of its growth  with $|R|$  is
stated (without proof) in Lemma 3 of \cite{Hoo}; an effective bound is:

\medskip\noindent\textbf{Lemma 6.} \textit{For any  integers $k,u,h_1,h_2$ with $h_1\le h_2$, we have}
\begin{equation}\left|\sum_{{h_1\le h\le h_2}}\!\!\!\!\!{}'\,\,e\left(\frac{uh'}{k}\right)
\right|< \left(\frac{k+h_2-h_1}{k}+2+2\ln(k)\right)\sqrt{k}\sqrt{\mathrm{gcd}(u,k)}\,d(k)\,,\label{lemma6}\end{equation}
\textit{where $h'$ denotes the inverse of $h$ (mod $k$) and the prime over the summation
means to restrict to {gcd}$(h,k)=1$.}

\medskip\noindent\textit{Proof.} Assume for now that $h_2- h_1< k$. 
For any integer $h$, $h_1\le h\le h_1+k-1$, the map sending $h$ to
\begin{equation}\left(\left(\frac{h-h_2-1/2}{k}
\right)\right)-\left(\left(\frac{h-h_1+1/2}{k}\right)\right)+\frac{h_2-h_1+1}{k} =\left\lfloor\frac{h-h_1+1/2}{k}
\right\rfloor-\left\lfloor\frac{h-h_2-1/2}{k}\right\rfloor\nonumber\end{equation}
equals 1 for $h_1\le h\le h_2$ and 0 otherwise. Now, we claim that for $0<x<1$, 
\begin{equation}\label{sawbound}
\left|\left(\left(x\right)\right)+\sum_{j=1}^k\frac{\sin(2\pi jx)}{\pi j}\right|< \frac{1}{2k\,\mathrm{min}(x,1-x)}\,.\end{equation}
To see this, it suffices to consider $0<x\le 1/2$. Write $K(x)=1+2\sum_{j=1}^k\cos(2\pi jx)=
{\sin(\pi\,(2k+1)x)}/{\sin(\pi x)}$; then 
\begin{align}\left(\left(x\right)\right)+&\sum_{j=1}^k\frac{\sin(2\pi jx)}{\pi j}=\int_0^xK(y)dy-1/2
=-\int_x^{1/2}K(y)dy\\=&\frac{-1}{\pi\,(2k+1)}\frac{\cos(\pi\,(2k+1)x)}{\sin(\pi x)}+
\frac{1}{2k+1}\int_x^{1/2}\frac{\cos(\pi\,(2k+1)y)}{\sin^2(\pi y)}dy\,.\end{align}
Since $\sin(\pi t)\ge 2t$ for $0\le t\le 1/2$, we get
\begin{equation}\left|\left(\left(x\right)\right)+\sum_{j=1}^k\frac{\sin(2\pi jx)}{\pi j}\right|\le
\frac{1}{2\pi\,(2k+1)x}+\frac{1}{4\,(2k+1)x}\end{equation}
which implies the weaker \eqref{sawbound}.

Write $S(v,u;k)$ for the (complete) Kloosterman sum $\sum'_{1\le h\le k}e(({vh+uh'})/{k})$
where the prime denotes restricting the sum to $\mathrm{gcd}(h,k)=1$.
Then e.g. Lemma 2 of \cite{Hoo2} gives an effective bound for it:
\begin{equation}\label{complKloos}
|S(v,u;k)|\le
\sqrt{k}\sqrt{\mathrm{gcd}(u,k)}\,d(k)\,.\end{equation}
Finally, it is elementary that $\sum_{j=1}^k1/j\le1+\ln(k)$.
Therefore, putting all this together, we obtain \begin{align}
\left|\sum_{h=h_1}^{h_2}\!{}'\,e\left(\frac{h'}{k}\right)\right|&\le
\frac{h_2-h_1+1}{k}|S(0,u;k)|+\left|\sum_{j=1}^k\frac{e(-(k+1/2)j/k)-
e(-j/(2k))}{2\pi j}S(j,u;k)\right|\nonumber\\ &\ \  +2\sum_{j=1}^{k/2}\frac{1}{2k\,(j-1/2)/k}
\le \frac{h_2-h_1+1}{k}\sqrt{k}\sqrt{\mathrm{gcd}(u,k)}\,d(k)\nonumber\\
&\ \ +\frac{1}{\pi}\sqrt{k}\sqrt{\mathrm{gcd}(u,k)}\,d(k)(1+\ln(k))+(1+\ln(k))\,.\end{align}
That the inequality \eqref{lemma6} also holds when $h_2-h_1\ge k$ now follows from
\eqref{complKloos}.\qquad
\textbf{QED to Lemma 6}\medskip

Effective bounds for the divisor function $d(n)=\sigma_0(n)$ are $d(n)\le C_\epsilon n^\epsilon$,
for any $\epsilon>0$, where 
\begin{equation}\label{dnbound}
C_\epsilon=\prod_{p<2^{1/\epsilon}}\frac{1}{\epsilon\,\ln(p)\,e^{1-\epsilon\ln(p)}}\,,\end{equation}
where the product runs over all primes $p<2^{1/\epsilon}$.
To see this, recall $d(n)=\prod_p(a_p+1)$ when $n=\prod_pp^{a_p}$ is the prime decomposition,
so $d(n)/n^\epsilon=\prod_p\frac{a_p+1}{p^{\epsilon a_p}}$. The primes appearing in
\eqref{dnbound} are precisely those for which $\frac{a_p+1}{p^{\epsilon a_p}}$ can be $>1$;
the power $a_p$ is then chosen to maximise this factor. Any $\epsilon<1/2$ works for us,
e.g. $C_{1/4}\approx 9.11...$ (in fact a slightly more refined analysis shows $C_{1/4}$ can be taken to
be 8.55). Also, from $\ln\,y\le y-1$ we obtain $\ln(k)<12\,k^{1/12}-12$. Therefore we obtain as a bound on \eqref{parsum}:
\begin{align}
&\left|\sum_{S_-/\delta\le S\le S_+/\delta} e\left(\frac{uS'}{R}\right)\varphi(r,\delta S)\right|
\nonumber\\&\ \le 
(103\left|\frac{rh}{\delta}\right|^{5/6}-74.8\left|\frac{rh}{\delta}\right|^{3/4})
( .0084\frac{n^5h^6}{r^2}\sqrt{|D|}+.000377\frac{n^8h^9}{|r|^3}\delta|D|+.0164\frac{n^{4}h^6}{\delta|r|^3}\sqrt{|D|} 
+\frac{.091n^3h^3}{r^2})
\nonumber\\&\ <(1.065\frac{n^5h^{41/6}}{\delta^{5/6}|r|^{7/6}}-.733\frac{n^5h^{27/4}}{\delta^{3/4}|r|^{5/4}})\sqrt{|D|}
+(.0389\frac{n^8h^{59/6}\delta^{1/6}}{|r|^{13/6}}-.028
\frac{n^8h^{39/4}\delta^{1/4}}{|r|^{9/4}})|D|
\,.\end{align}
We can now bound the limit of \eqref{sum11} as $X\rightarrow\infty$, using the inequalities
$n^{1-\mu}/(\mu-1)<\sum_{r=n}^\infty r^{-\mu}<n^{-\mu}+n^{1-\mu}/(\mu-1)$ for $\mu>1$ which follow
by comparing series with integral:
\begin{align}&\mathrm{lim}_{X\rightarrow\infty}\left|\sum_>\right|\le  (19.136n^{29/6}h^{41/6}
-9.67n^{19/4}h^{27/4})\sqrt{|D|}
+(.455n^{41/6}h^{59/6}-.359n^{27/4}h^{39/4})|D|\,.
\end{align}
When $n=2$ (so $h=\delta=1$), we get instead the weaker lim$_{X\rightarrow\infty}|\sum_>|<
2136\sqrt{|D|}+61|D|$.


Bounding $\sum_<$ is completely analogous; \eqref{separ} becomes
\begin{equation}
\frac{m}{2nc}\equiv_1\frac{ur'}{S}-\frac{\alpha r+\beta s}{2ns\,(\alpha r^2+\beta rs+\gamma n's^2)}\,,\end{equation}
where $\delta=\mathrm{gcd}(r,h)\in\{1,3\}$, $S=ns/\delta$, $r'$ is the inverse of $r/\delta$
mod $S$, and $u$ is defined mod $S$ by $u\equiv_S\delta^{-1}\,(1+n's^2-n'{}^2s^4)$.
The other equations and bounds are obtained by interchanging $r\leftrightarrow s$ and
$\alpha\leftrightarrow n'\gamma$. In \eqref{sum11} and elsewhere it suffices to take $\delta\in\{1,3\}$. Then \eqref{deltaphi},\eqref{phi} become
\begin{align}
&|\varphi(r,s)-\varphi(r+\delta,s)|\le .056\frac{n^2h^6\delta}{|s|^{3}}\sqrt{|D|}+.0093\frac{n^{3}h^{9}\delta^2}{s^{4}}|D| +.34\frac{nh^6}{s^4}\sqrt{|D|}\,,\\
&\ \ \ \ \ \ \ \ \ |\varphi(\delta\lfloor S_+/\delta\rfloor,s)|\le nh^3s^{-2}/3\,.\end{align}
Therefore \eqref{parsum} becomes
\begin{align}
&\left|\sum_{S_-/\delta\le R\le S_+/\delta} e\left(\frac{uR'}{S}\right)\varphi(\delta R,s)\right|
<(1.916\frac{n^{29/6}\delta^{1/6}}{|s|^{13/6}}-1.26\frac{n^{19/4}\delta^{1/4}}{|s|^{9/4}})h^9|D|
\nonumber\\&\qquad\qquad\qquad+( 12.98\frac{n^{23/6}h^6}{\delta^{5/6}|s|^{7/6}}+70.1\frac{n^{17/6}h^6}{\delta^{11/6}|s|^{13/6}}
-8.56\frac{n^{15/4}h^6}{\delta^{3/4}|s|^{5/4}}-46\frac{n^{11/4}h^6}{\delta^{7/4}|s|^{9/4}})\sqrt{|D|}
\,,\nonumber\end{align}
and we obtain
\begin{align}&\mathrm{lim}_{X\rightarrow\infty}\left|\sum_<\right|\le(7.832n^{29/6}h^9-2.33n^{19/4}h^9)|D|
+(201n^{23/6}h^6-72n^{15/4}h^6)\sqrt{|D|}\,.\nonumber
\end{align}

All that remains is the sum over the representatives $[n\alpha,\beta,\gamma/h]$ of
$\Gamma_0(n;h)$-equivalence classes of quadratic forms of discriminant $D$. Let $h'_{n;h}(D)$ denote the
number of those (not necessarily primitive) equivalence classes. Note that our sum is now independent
of $\alpha,\beta,\gamma$, so a bound for $|Z_{n/2;h}(3/4)|$ will be the sum of our bounds for 
$\sum_>$ and $\sum_<$, multiplied by $h'_{n;h}(D)/2$ (the 2 here compensates for the overcounting,
as mentioned earlier). Now, since all stabilisers here are $\pm I$, $h'_{n;h}(D)$ will equal the
number of $\Gamma_0(n')$-equivalence classes of forms $[n'\alpha,\beta,\gamma]$, times
the index of $\Gamma_0(n;h)$ in $\Gamma_0(n|h)$, and so we obtain the (crude) bound
\begin{equation}
h'_{n;h}(D)\le h'(D)\,\left|\frac{\mathrm{SL}_2(\bbZ)}{\Gamma(n,h)}\right|\label{hpbound1}\end{equation}
where $h'(D)$ is the total number of SL$_2(\bbZ)$-equivalence classes of integral forms
(not necessarily primitive) of discriminant $D$, and where $\Gamma(n,h)$ consists of all
matrices $\left({a\atop c}{b\atop d}\right)$ where $n$ divides $c$ and $h$ divides $b$
(this clearly forms a group). In deriving \eqref{hpbound1}, we are conjugating everything by
$\left({h\atop 0}{0\atop 1}\right)$; $\Gamma(n,h)$ arises as the conjugate of the subgroup
$\Gamma_0(nh)$ of $\Gamma_0(n;h)$. Now, the index of $\Gamma(n,h)$ in SL$_2(\bbZ)$
equals that of its conjugate $\Gamma_0(nh)$, since indices of Fuchsian groups can be expressed as a ratio of
areas of fundamental domains and conjugating by SL$_2(\bbR)$ preserves those areas.
Let $h(D)$ denote the class number of $D$, i.e. the number of SL$_2(\bbZ)$-equivalence
classes of \textit{primitive} integral quadratic forms of discriminant $D$. Then we
have the bound 
\begin{equation} h(D)< \frac{\sqrt{|D|}}{\pi}\,(2+\log|D|)\,\end{equation}
(see Theorems 12.14.3 and 12.10.1 in \cite{Hua}). Hence  
\begin{align}
\frac{h'(D)}{2}&=\frac{1}{2}\sum_{m^2|D}h(D/m^2)<\frac{\sqrt{|D|}}{2\pi}\sum_{m^2|D}\frac{2-\log\,m}{m}+\sqrt{|D|}\frac{\log{|D|}}{4\pi}
\sum_{m^2|D}1\nonumber\\&<0.32\sqrt{|D|}+.681|D|^{5/8}\log\,|D|\,,\end{align}
using our earlier bound for the divisor function $d(n)$. Putting this all together, we obtain
\begin{equation}h'_{n;h}(D)/2<(0.32\sqrt{|D|}+.681|D|^{5/8}\log\,|D|)nh\prod_{p|(nh)}\frac{p+1}{p}\,,\end{equation}
and hence the bound given in the statement of our theorem. \textbf{QED to Theorem 3}\medskip

We can approximate the $2.13|D|^{1/8}\log|D|$ with the upper bound $120(|D|^{1/7}-|D|^{1/8})$
by the usual reasoning.
Putting this all together, we obtain the following (very crude) bounds for the character values
for each $g\in M_{24}$ and $k\ge 150$:
\begin{align}H_k(1)&>\frac{4}{K}e^{\pi K/2}-\frac{4\pi}{\sqrt{2}}e^{\pi K/4}-2.5\times 10^4K^{23/7}\,,\label{boundg1}\\  \label{boundgne1}
|H_k(g)|&<\frac{4}{K\sqrt{n_g}}e^{\pi K/(2n_g)}+{\sqrt{8}\pi}\,e^{\pi K/(4n_g)}+a_g\times 10^{b_g}K^{23/7}\,,\end{align}
where $n_g=|g|$ is the order of $g$, $K=\sqrt{8k-1}$, and the values $a_g,b_g$ are computed
from Theorem 3 and are collected in Table 2. We find that when $k\ge 390$, the ratio
$(H_k(1)/|M_{24}|)/(|H_k(g)|/|C_{M24}(g)|)$ exceeds $1.6\times 10^5$ for all elements $g\ne 1$ except for
$g$ in conjugacy class 12B, where the ratio is $>1.3$. In particular we find that 
\begin{equation}\frac{H_k(1)}{|M_{24}|}>\sum_{g\ne 1}\frac{|H_k(g)|}{|C_{M24}(g)|}\,,\label{thebound}\end{equation}
for all $k\ge 390$, where the sum is over all conjugacy classes in $M_{24}$.
As explained in \eqref{posbound}, this is sufficient to deduce positivity of the multiplicities
of all $M_{24}$-irreps in these class functions $H_k$. 
Again, positivity (or rather nonnegativity) has been obtained experimentally for $n$ up to 500, far more than we need.

It is somewhat disappointing that we need to check by hand (or rather, by computer!)
positivity for so many $k$, when empirically \eqref{thebound} is satisfied for $k>30$.
 The cause of this are the bounds of Theorem 3, which apparently are far looser than they could
be. But the value gained in tightening these bounds is far less than the effort to be spent,
since it is so easy to calculate these character multiplicities. The reader interested in reducing 
390 should focus on improving the bound for the case 12B (i.e. fix $n=24,h=12$ from the
start, and restrict to $k>100$ say; the index of $\Gamma_0(nh)$ in $\Gamma_0(n;h)$ here is 12, so the Theorem 3
bound can be improved immediately by that factor 12).

The bounds \eqref{boundg1},\eqref{boundgne1} imply:

\medskip\noindent\textbf{Theorem 4.} \textit{Choose any $M_{24}$-irrep $\rho$. Then the
multiplicity mult$_\rho(H_n)$ of $\rho$ in $H_n$ tends to $\infty$ as $n\rightarrow\infty$.
Moreover, mult$_\rho(H_n)>0$ for all $n\ge 25$.}\medskip

Indeed, in this section we proved positivity of the multiplicities for $n\ge 390$; the values $25\le n<390$
can be checked explicitly. An immediate corollary of Theorem 4 is the validity of Conjecture 5.11
in Umbral Moonshine \cite{CDH}. However it should be remarked that the proof of this in
\cite{CHM} actually established a much stronger statement. Oddness of McKay-Thompson
coefficients is far less trivial than strict positivity of certain  multiplicities, since as we see
almost every multiplicity will be strictly positive.

\section{Is the Conway group the stringy symmetry?}

In beautiful work, \cite{GHV3} followed Kondo's treatment \cite{Ko} of Mukai's classification
\cite{Mu} of symplectic automorphisms of K3 surfaces, to obtain the symmetries of
K3 sigma models. (Similar considerations are considered in \cite{TW}.) It turns out all these stringy symmetries are subgroups of the
Conway groups Co$_0$ or Co$_1$ (but not $M_{24}$), where they will necessarily generate
the full group Co$_0$ (resp. Co$_1$). This begs the (perhaps naive) question: should the
automorphism group of Mathieu Moonshine actually be a Conway group?

The Conway group Co$_0$ is the automorphism group of the Leech lattice; the simple group 
Co$_1$ is its quotient by its centre $\pm 1$. Hence the Co$_1$-irreps are a subset of
those of Co$_0$, consisting of those Co$_0$-irreps whose kernel contains $-1$. We call
a Co$_0$-irrep a \textit{spinor} if it is not trivial on $-1$.   
$M_{24}$ is a subgroup of both Co$_0$ and Co$_1$. In this section we give the 
 restrictions of Co$_0$ (hence Co$_1$) irreps to $M_{24}$, for irreps with dimension up to 1 million 
(there are  32 of these).  The notation comes from the Atlas 
of Finite Groups \cite{Atl}; in a couple places it differs slightly from 
that of  \cite{EOT}
(e.g. the dimension-1035 irreps are in a different order: the conjugate ones we write as $1035'$ and $1035''$). 
 The first step was identifying the conjugacy classes; once these are known,  the 
character multiplicities follow directly from 
the chararacter table of $M_{24}$ together with the orthogonality relations. The arguments are
straightforward and we'll avoid the details, giving only the results.

\bigskip
\centerline {{\bf Table 5.} Matching of conjugacy classes}
$${\scriptsize\vbox{\tabskip=0pt\offinterlineskip
  \def\tablerule{\noalign{\hrule}}
  \halign to 7in{
    \strut#&
    \hfil#&\vrule#&\vrule#&\hfil#&\vrule#&    
\hfil#&\vrule#&\hfil#&\vrule#&\hfil#&\vrule#&\hfil#&\vrule#&\hfil#&\vrule#&\hfil#&\vrule#&\hfil#&\vrule#&\hfil#&\vrule#&\hfil#&\vrule#&\hfil#&\vrule#&\hfil#&\vrule#&\hfil#&\vrule#&\hfil#&\vrule#&\hfil#&\vrule#&\hfil#&\vrule#&\hfil#&\vrule#&\hfil#&\vrule#&\hfil#&\vrule#&\hfil#&\vrule#&\hfil#&
\vrule#&\hfil#&\vrule#&\hfil#&\vrule#&\hfil#&\vrule#&\hfil#&\vrule#&\hfil#
\tabskip=0pt\cr
&$M_{24}\,\,$&\,&&\,\,1A\,\,&&\,\,2A\,\,&&\,\,2B\,\,&&\,\,3A\,\,&&\,\,3B\,\,&&\,\,4A\,\,&&\,\,4B\,\,&&\,\,4C\,\,&&\,\,5A\,\,&&\,\,6A\,\,&&\,\,6B\,\,&&\,\,7A\,\,&&\,\,7B\,\,&&\,\,8A\,\,&&\,\,10A\,\,&&\,\,11A\,\,&&\,\,12A\,\,&&\,\,12B\,\,&&\,\,14A\,\,&&\,\,14B\,\,&&\,\,15A\,\,&&\,\,15B\,\,&&\,\,21A\,\,&&\,\,21B\,\,&&\,\,23A\,\,&&\,\,23B\,\,\cr
\tablerule&Co$_1\,\,$&\,&&\,\,1A\,\,&&\,\,2A\,\,&&\,\,2C\,\,&&\,\,3B\,\,&&\,\,3D\,\,&&\,\,4D\,\,&&\,\,4C\,\,&&\,\,4F\,\,&&\,\,5B\,\,&&\,\,6E\,\,&&\,\,6I\,\,&&\,\,7B\,\,&&\,\,7B\,\,&&\,\,8E\,\,&&\,\,10F\,\,&&\,\,11A\,\,&&\,\,12J\,\,&&\,\,12M\,\,&&\,\,14A\,\,&&\,\,14B\,\,&&\,\,15D\,\,&&\,\,15D\,\,&&\,\,21C\,\,&&\,\,21C\,\,&&\,23A\,\,&&\,\,23B\,\,\cr\noalign{}}}}$$

\medskip
\centerline {{\bf Table 6.} The 32 smallest Conway irreps}
\begin{align}&&{\scriptsize\vbox{\tabskip=0pt\offinterlineskip
  \def\tablerule{\noalign{\hrule}}
  \halign to 6.18in{
    \strut#&
    \hfil#&\vrule#&\hfil#&\vrule#&    
\hfil#&\vrule#&\hfil#&\vrule#&\hfil#&\vrule#&\hfil#&\vrule#&\hfil#&\vrule#&\hfil#&\vrule#&\hfil#&\vrule#&\hfil#&\vrule#&\hfil#&\vrule#&\hfil#&\vrule#&\hfil#&\vrule#&\hfil#&\vrule#&\hfil#&\vrule#&\hfil#&\vrule#&\hfil#&\vrule#&\hfil#
\tabskip=0pt\cr
&$\tsi_1\,\,$&&$\,\,\si_1\,\,$&&\,\,$\si_2$\,\,&&\,\,$\si_3$\ \,&&\,\,$\tsi_2$\ \,&&\,\,$\tsi_3$\ \,&&\,\,$\tsi_4$\ \,&&\,\,$\si_4$\ \,\,&&\,\,$\si_5$\ \,\,&&\,\,$\si_6$\ \,\,&&\,\,$\si_7$\ \,\,&&\,\,$\tsi_5$\ \,\,&&\,\,$\si_8$\ \,\,&&\,\,$\si_9$\ \,\,&&\,\,$\si_{10}$\ \,\,&&\,\,$\tsi_6$\ \,\,&&\,\,$\tsi_7$\ \,\,&&\,\,$\tsi_8\ \,\,$\cr
\tablerule&$24\,\,$&&$\,\,276\,\,$&&\,\,299\,\,&&\,\,1771\,\,&&\,\,2024\,\,&&\,\,2576\,\,&&\,\,4576\,\,&&\,\,8855\,\,&&\,\,17250\,\,&&\,\,27300\,\,&&\,\,37674\,\,&&\,\,40480\,\,&&\,\,44275\,\,&&\,\,80730\,\,&&\,\,94875\,\,&&\,\,95680\,\,&&\,\,170016\,\,&&\,\,299000\cr
\noalign{}}}}\nonumber\\
&&{\scriptsize\vbox{\tabskip=0pt\offinterlineskip
  \def\tablerule{\noalign{\hrule}}
  \halign to 5.6in{
    \strut#&   
    \hfil#&\vrule#&\hfil#&\vrule#&    
\hfil#&\vrule#&\hfil#&\vrule#&\hfil#&\vrule#&\hfil#&\vrule#&\hfil#&\vrule#&\hfil#&\vrule#&\hfil#&\vrule#&\hfil#&\vrule#&\hfil#&\vrule#&\hfil#&\vrule#&\hfil#\tabskip=0pt\cr     
&$\si_{11}$\ \,&&\,\,$\tsi_9$\ \,\,&&\,\,$\si_{12}$\ \,&&\,\,$\tsi_{10}$\ \,&&\,\,$\si_{13}$\ \,&&\,\,$\tsi_{11}$\ \,\,&&\,\,$\si_{14}$\ \,\,&&\,\,$\si_{15}$\ \,\,&&\,\,$\si_{16}$\ \,\,&&\,\,$\si_{17}$\ \,\,&&\,\,$\tsi_{12}$\ \,\,&&\,\,\,$\si_{18}$\ \,\,&&\,\,$\si_{19}\ \,\,$\cr
\tablerule&313950\,\,&&\,\,315744\,\,&&\,\,345345\,\,&&\,\,351624\,\,&&\,\,376740\,\,&&\,\,388080\,\,&&\,\,483000\,\,&&\,\,644644\,\,&&\,\,673750\,\,&&\,\,673750\,\,&&\,\,789360\,\,&&\,\,822250\,\,&&\,\,871884\cr
\noalign{}}}}\nonumber\end{align}

We find that any  \textit{virtual character} $\rho=m_0+m_1\rho_1+\cdots+
m_{20}\rho_{20}$ of $M_{24}$, i.e. any linear combination with integer coefficients $m_i$ 
of $M_{24}$ irreps, is the restriction of a virtual Co$_0$-representation 
involving  the first 32 irreps of Co$_0$, iff the multiplicities obey $m_2=m_{\overline{2}}$, $m_3=m_{\bar{3}}$, 
$m_8=m_{\overline{8}}$, $m_{10}=m_{\overline{10}}$, as well as the relation 
$m_9+m_{14}\equiv_2m_{15}$. Thus it is hard for an $M_{24}$-representation to be
a restriction of a virtual Co$_0$-representation (0\% 
will be). Although the restriction map is not
surjective, Theorem B tells us the Mathieu-Moonshine characters $H_n$  lie in the image
of that restriction. This is the content of Theorem C given earlier. In particular we have:
{\small\begin{eqnarray}
&&\rho_1=\tsi_1 - 1\,,\nonumber\\
&&\rho_2+\overline{\rho_2}=  37 \si_{14} + 39 \si_{13} + 28 \tsi_{10 }+ \si_{12}+ 12 \si_{11}  + 23 \tsi_8 + \si_{10}+ 25 \si_9 + 2\si_5    + 194 \tsi_4 + 24 \tsi_3+ 18 \si_3 \nonumber\\&&\qquad + 269 \tsi_1- 28 \si_{18} - 11 \tsi_{12}  - 42 \tsi_9- 44 \tsi_6 - 9 \si_8  - 100 \tsi_5- 18 \si_7  - 78 \si_6 - 144 \si_2 - 218 \si_1- 155
\,,\nonumber\\
&&\rho_3+\overline{\rho_3}=  14\si_{13} + 11\tsi_{10}+ 15\si_{12}    +4\si_{11}+ 6\tsi_8 + 15 \si_{10} + 8\si_9  + 105\tsi_4 + 14\tsi_3+ 64 \si_3 + 121\tsi_1 \nonumber\\ &&\qquad
- 10\si_{18}- 4\tsi_{12}- 15\tsi_9 - 15\tsi_6 - 20 \si_8- 12\tsi_5- 20\si_7 - 13\si_6 - 97\si_2 - 49\si_1 - 27 \,,\nonumber\\
&&\rho_4= \si_2+1-2\tsi_1\,,\nonumber\\
&&\rho_5= \si_1+ 1 -\tsi_1 \,,\nonumber\\
&&\rho_6=   8 \si_{14}   + 14 \si_{13}+ 11 \tsi_{10} + 7 \si_{12}+ 4 \si_{11}+ 6 \tsi_8 
 + 7 \si_{10}   + 8 \si_9+ 89\tsi_4  + 14 \tsi_3+32\si_3 + 121\tsi_1\nonumber\\ &&\qquad
 - 10 \si_{18} -4 \tsi_{12} - 15 \tsi_9   - 15\tsi_6
 - 12 \si_8- 28 \tsi_5- 12 \si_7- 21\si_6    - 81 \si_2- 73 \si_1- 51  \,,\nonumber\\
&&\rho_7=36\si_{14}   + 41\si_{13}+ 33 \tsi_{10} + 6\si_{12} + 12\si_{11}
 + 25\tsi_8  + 6\si_{10}+ 25\si_9 + 258\tsi_4+ 50\tsi_3 +55\si_3  + 353 \tsi_1 
   \nonumber\\ &&\qquad - 30\si_{18}- 12\tsi_{12} - \si_{17}   
  - 43 \tsi_9 - 45\tsi_6 - 29\si_8- 100\tsi_5 - 32\si_7 - 86\si_6 - 2\si_5- 226\si_2- 256 \si_1- 182 \, ,\nonumber\\
&&\overline{\rho_7}= \si_{17}+ \si_{14}   + 41\si_{13} + 39\si_{12}  + 25 \tsi_{9}  + 12\si_{11} + 23 \tsi_8 
 + 39 \si_{10} + 27\si_9+ 2 \si_5 + 307 \tsi_4 + 40 \tsi_3+181\si_3\nonumber\\ &&\qquad+ 343\tsi_1- 28\si_{18}- 12 \tsi_{12}- 41 \tsi_9 - 47 \tsi_6- 52 \si_8  - 31\tsi_5
  - 61\si_7   - 49 \si_6 - 278 \si_2- 155\si_1-85\,,\nonumber\\
&&\rho_8+\overline{\rho_8}=5 \si_{14} + 6 \si_{13}+ 4\tsi_{10}  + 2\si_{12} + 2\si_{11} + 2 \tsi_8 
     + 3\si_{10} + 4\si_9 +\si_5  + 44 \tsi_4+ 7\tsi_3+12\si_3 + 66\tsi_1  \nonumber\\ &&\qquad
 - 4\si_8 - 4\si_{18}- 2\tsi_{12} - 6\tsi_9 - 7\tsi_6 - 16\tsi_5- 6\si_7 - 13\si_6- 44\si_2- 46\si_1  - 33 \,,\nonumber\\
&&\rho_9+\rho_{15}= 12\si_{18}+ 5\tsi_{12}+ 17\tsi_9+ 20\tsi_6  + 22\si_8 + 18\tsi_5+ 25\si_7+ 25\si_6  + 125 \si_2 + 82\si_1 + 51
 \nonumber\\ &&\qquad   - 3\si_{14}  - 17\si_{13}- 12\tsi_{10} - 14\si_{12}
   - 5\si_{11}- 10\tsi_8  - 14\si_{10}- 11\si_9- 133\tsi_4- 23\tsi_3-69 \si_3- 159 \tsi_1\,,\nonumber\\
&&\rho_{10}+\overline{\rho_{10}}=22\si_{18}+ 10 \tsi_{12}+ 31\tsi_9+ 36\tsi_6 + 53\si_8+ 24\tsi_5 + 55\si_7+ 44\si_6
          +\si_5
+ 270\si_2+ 150\si_1+ 87\nonumber\\ &&\qquad- 32\si_{13}- 22\tsi_{10} - 33\si_{12}- 10\si_{11}  - 17\tsi_8- 34\si_{10}- 20\si_9- 278 
\tsi_4 - 50\tsi_3-161\si_3- 321\tsi_1 \,,\nonumber\\
&&\rho_{11}= 12\si_{18}+ 5\tsi_{12}+ 17\tsi_9 + 19\tsi_6  + 19 \si_8+ 28\tsi_5+ 20\si_7+ \si_5 + 30\si_6
     + 116\si_2+ 94\si_1+ 63\nonumber\\ &&\qquad- 8\si_{14} - 17\si_{13}- 12\tsi_{10}- 9\si_{12}     - 10\tsi_8- 5\si_{11}- 9\si_{10}- 11\si_9
- 125\tsi_4 - 23\tsi_3      -51\si_3- 158 \tsi_1 \,,\nonumber\\
&&\rho_{12}=\si_3\,,\nonumber\\
&&\rho_{13}=14\si_{13}+ 10\tsi_{10}  + 14\si_{12}+ 4 \si_{11}
     + 8\tsi_8+ 14\si_{10}+ 9\si_9+ 110\tsi_4+ 18 \tsi_3+65 \si_3+ 124\tsi_1\nonumber\\ &&\qquad- 10\si_{18}
     - 4\tsi_{12}- 14\tsi_9 - 16\tsi_6 - 21\si_8- 10\tsi_5- 22\si_7 - 16\si_6-\si_5 - 102\si_2- 53\si_1- 30 \,,\nonumber\\
&&\rho_{14}+\rho_{15}=  12\si_{18}+ 5\tsi_{12} + 18\tsi_9 + 19\tsi_6 + 21\si_8+ 16\tsi_5+ 23\si_7 + 21\si_6+ 111\si_2+ 65\si_1+ 38\nonumber\\ &&\qquad
- 2\si_{14}- 17\si_{13}- 12\tsi_{10}- 15\si_{12}- 5\si_{11}- 10 \tsi_8 - 15\si_{10} - 11\si_9- 124 \tsi_4
- 17\tsi_3 -69\si_3- 141\tsi_1  
\,,\nonumber 
\end{eqnarray}}
{\small\begin{eqnarray}
&&2\rho_{15}=  30\si_{14}+ 39\si_{13}+ 28\tsi_{10}+ 11\si_{12}+ 12\si_{11}+ 22\tsi_8+ 11\si_{10}+ 25 \si_9 + 244\tsi_4+ 40\tsi_3 + 67\si_3\nonumber\\ &&\qquad + 326\tsi_1
 - 28\si_{18} - 12\tsi_{12}- 40\tsi_9 - 44\tsi_6- 26\si_8- 90\tsi_5-33\si_7- 76\si_6   - 211\si_2- 229\si_1- 158
\,,\nonumber \\
&&\rho_{16}=\tsi_4+4\tsi_1 -2\si_2-2\si_1-2\,,\nonumber\\
&&\rho_{17}=5\si_{18} + 2\tsi_{12} + 7\tsi_9 + 7\tsi_6  + 12\si_8+ 5\tsi_5+11\si_7+ 8\si_6 +\si_5 
   + 54\si_2+ 27\si_1 + 15\nonumber\\ &&\qquad- 7\si_{13}- 5\tsi_{10}- 7\si_{12}- 2\si_{11}- 4\tsi_8- 7\si_{10}- 4\si_9
     - 57\tsi_4- 10\tsi_3 - 34\si_3 - 64\tsi_1\,,\nonumber\\
&&\rho_{18}=\si_{14}+ 6\si_{13}+ 4\tsi_{10}+ 6\si_{12}+ 2\si_{11}+ 2\tsi_8+ 6\si_{10}+ 4\si_9 
     + 46\tsi_4+ 7\tsi_3  + 24\si_3+ 58\tsi_1\nonumber\\ &&\qquad - 4\si_{18} - 2\tsi_{12}- 6\tsi_9 - 7\tsi_6- 7\si_8
     - 8\tsi_5-8\si_7- 8\si_6- 44\si_2- 30\si_1- 19\,,\nonumber\\
&&\rho_{19}=8\si_{14}+ 10\si_{13}+ 7\tsi_{10}+ 2\si_{12}+ 3\si_{11}+ 6\tsi_8 + 2\si_{10}+ 6\si_9 
     + 68\tsi_4+ 14\tsi_3 + 18\si_3+ 95\tsi_1\nonumber\\ &&\qquad- 7\si_{18}- 3\tsi_{12}    - 10\tsi_9- 11\tsi_6- 8\si_8- 23\tsi_5-9\si_7 - 21\si_6
  - 64\si_2- 67\si_1- 48\,,\nonumber\\
&&\rho_{20}= 41\si_{13}+ 29\tsi_{10}+ 41\si_{12}+ 12\si_{11}+ 24\tsi_8+ 41\si_{10} + 26\si_9+ 318\tsi_4   + 45\tsi_3 + 191\si_3+ 346\tsi_1\nonumber\\ &&\qquad
   - 29\si_{18} - 12\tsi_{12}- 42\tsi_9- 46\tsi_6 - 59\si_8- 29\tsi_5-64\si_7- 49\si_6
 - 287\si_2- 149\si_1- 78\,.\nonumber
\end{eqnarray}}
We considered more Co$_0$-irreps than $M_{24}$ ones; the additional relations between restrictions are:
{\small\begin{eqnarray}
&&\tsi_2 +\tsi_1 =\si_3 +\si_1+1      \,,\nonumber\\
&&\tsi_5+\tsi_4+\si_3+1=\si_7+\si_4+\si_2 \,,\nonumber\\
&& \tsi_{11} +\tsi_3=\tsi_8 +\tsi_5   +\si_7+\si_4+ 2\tsi_2+\si_2+\si_1+\tsi_1\,,\nonumber\\
&& \si_{14}+\si_8=\tsi_8 +\tsi_7 +\si_7 + 2\si_4+\tsi_2+\si_2+ 2\si_1\,,\nonumber\\
&& \tsi_8 +\tsi_7+\si_2=\si_{12} +\si_{10}+\si_6+\si_3+\tsi_1  \,,\nonumber\\
&&  \si_{19} +\si_8= \si_{18}+ 2\si_7+ 2\si_4+\si_2+2\si_1\,,\nonumber\\
&&\tsi_6 +4\tsi_4+ 3\tsi_3 + 24535 
\tsi_2   + 24541\tsi_1=   2\si_8+2\si_5 + 24533\si_3+6\si_2+ 24537\si_1+ 24537\,,\nonumber\\
&& 85\si_{14}+ 88\si_{13}+ 62\tsi_{10}+ 4\si_{12}+ 26\si_{11}+ 50\tsi_8+ 4\si_{10}+ 56\si_9 + 516\tsi_4+ 97\tsi_3 + 72\si_3 + 754\tsi_1 \nonumber\\ &&\qquad=62\si_{18}+ 26\tsi_{12}+90\tsi_9
+ 99\tsi_6+ 41\si_8 + 234\tsi_5+53\si_7 + 190\si_6  + 453\si_2+ 582\si_1+ 428\,.\nonumber\end{eqnarray}}

Of course we'd prefer true representations to virtual ones, but perhaps we should not be too
surprised that virtual representations arise here, because the definition \eqref{ellge}
of elliptic genus involves
signs, and after all even $H_{00}$ and $H_0$ were virtual (but see Section 6).

However, the dimensions of the virtual Co$_0$-representations needed will be \textit{extremely} 
large. 
By the \textit{total dimension} of a virtual representation $\rho_+\ominus\rho_-$ we
mean the quantity dim$(\rho_+)+\mathrm{dim}(\rho_-)$.
 For any $M_{24}$ representation,
if there is any Co$_0$-virtual
representation which restricts to it, there will be infinitely many (just add to it any of the above
relations). For each basic combination
of $M_{24}$-irreps, we have selected above the Co$_0$-virtual representation of smallest total
dimension we could find, which restricts to it. (This ambiguity would be eliminated by
identifying if possible the twining elliptic genera $\phi_g$ for all $g\in \mathrm{Co}_0$.)

For example,  of all Co$_0$-irreps of dimension less than 1 million, only 2 of them contain $\rho_2$ 
or $\overline{\rho_2}$ as a summand: namely, 
a Co$_1$-irrep of dim 313950, and a spinor irrep of dim 789360. 
The smallest virtual Co$_0$-representation we could find which restricts to $\rho_2+\overline{\rho_2}$ (a 90-dimensional representation of $M_{24}$)
has total dimension over 100 billion.

In fact more is true. Suppose for contradiction the natural identification of the Witten index
$w_g$ with $\widetilde{\sigma}_1$ (the only Co$_0$-rep which restricts to $1+\rho_1$).
Now,  \cite{GHV3} compute (among other things) the
 twining characters for  the automorphism of the Gepner model $(1)^6$; they found two symmetries of order 9 (the bottom two rows of their Table 1), both of whose Witten indices
 equal 3, and whose twining elliptic genera are nonetheless different.
 However, Co$_0$ has only one conjugacy class
(namely class 9C) of order 9 with $\widetilde{\sigma}_1=3$, and this would imply those
twining genera should be equal. (We thank Roberto Volpato 
for sharing this observation.) This contradiction means that $w_g$ cannot equal $\widetilde{\sigma}_1(g)$, and so the Witten genus too will be virtual, having total dimension
much greater than 24. These large total dimensions 
 make the proposed extension to the Conway groups seem unlikely.

\section{Speculations}

Because the phenomena underlying Mathieu Moonshine are still obscure, we conclude
 with assorted questions and speculations. The most obvious challenge suggested by
this paper is to explicitly construct these representations $H_n$. The evidence that there is some 
vertex algebra-like object underlying the Mathieu Moonshine observations now seems overwhelming;
perhaps the most satisfying way to construct the $H_n$ would be to construct this vertex (super)algebra.
This would provide the algebraic underpinning of Mathieu Moonshine, and it should hint at
 its still-mysterious geometric and physical meanings.
 
 We regard this vertex superalgebra construction as the most important challenge of Mathieu Moonshine.
At $c=24$ we have an $N=0$ VOA (namely the Moonshine Module $V^\natural$) 
with lots of nice properties including an action of $\mathbb{M}$. At $c=12$ we have an
$N=1$ VOSA (namely Duncan's algebra \cite{Du}), with lots of nice properties including
an action of the Conway group. Could there be at $c=6$ an $N=2$ or $N=4$ VOSA
(namely the algebra underlying our Mathieu Moonshine) with lots of nice properties including
an action of $M_{24}$? After all, $M_{24}$, Co$_1$, $\mathbb{M}$ is the Holy Trinity of
sporadic finite simple groups (e.g. Griess \cite{Gr} constructed the Monster by starting with $M_{24}$,
lifting to Conway and then moving on to the Monster). (I thank Gerald H\"ohn and Chongying
Dong for informal discussions on this point).

Another possible construction, which may also bring 
 K3 into Mathieu Moonshine, is through the chiral 
de Rham complex \cite{MSV} (a sheaf of vertex superalgebras) associated to a K3 surface.  The trace over its global section (which is itself a vertex superalgebra) recovers elliptic genus
\cite{BoLi}. The orbifold theory, including the construction of twisted modules, 
has also been studied \cite{FrSz}.

One of the most intriguing aspects of Mathieu Moonshine is the \textit{positivity}
 of the irrep multiplicities. This discussion is already anticipated in Chapter 7.2 of Wendland's thesis \cite{We}. There we find that the partition function in the Ramond-Ramond sector of
 an $N=4$ $c=6$ superconformal field theory is
 \begin{align}\label{partfunc}
 Z_{RR}=&|ch^s_{0,0}|^2+h|ch_{1/4,0}^s|^2+\epsilon\,(ch^s_{0,0}ch_{1/4,0}^{s*}+c.c.)
 +Fch^l_{1/4,1/2}ch^{s*}_{0,0}+F'ch^s_{0,0}ch^{l*}_{1/4,1/2}\nonumber\\ &+Gch^l_{1/4,1/2}ch^{s*}_{1/4,0}+G'ch^s_{1/4,0}ch^{l*}_{1/4,1/2}+H|ch^l_{1/4,1/2}|^2\,,\end{align}
 using the same conventions for $N=4$ superconformal characters as \eqref{ellgenN4}, where $h,\epsilon$ are nonnegative integers
 and $F,F',G,G',H$ are functions of $\tau,z$ with only nonnegative integer coefficients
 in $q,\overline{q}$. The elliptic genus is
 \begin{equation}
 \phi=(\epsilon-2)ch^s_{0,0}+(h-2\epsilon)ch^s_{1/4,0}+(G-2F)ch^l_{1/4,1/2}\,.\end{equation}
For the K3 component of moduli space, $\epsilon=0$ and $h=h^{1,1}=20$.
In particular note
that the function $F$ counts holomorphic fields which exist in that theory but which are
not  contained in the $N=4$ vacuum representation. We can think of holomorphic fields, very roughly, as
additional symmetries of the theory; since we would expect that generically the chiral algebra
of the K3 sigma model is simply $N=4$ superconformal,
we should have generically $F=0$. This would mean the non-BPS part of the elliptic
genus is (generically) this function $G$  which has only non-negative coefficients.
From this the conjectured non-negativity would follow. This argument also applies to the
character-valued elliptic genus, so all multiplicities should be nonnegative. A similar
argument can be found in Ooguri \cite{Oo}. (We thank Katrin Wendland for discussions
on this point.)

It is commonly expected that the moduli space of $c=6$ $N=4$ superconformal field theories
 consists of two components: a toroidal one with vanishing elliptic
genus, and the K3 sigma models with elliptic genus $2\phi_{0,1}$.  But if this is true, then
there can be no $c=6$ $N=4$ superconformal field theory underlying Mathieu moonshine:
the work of \cite{GHV3} shows that the corresponding automorphism groups are too small.

But  $N=4$ (more precisely, $N=(4,4)$) theories possess more supersymmetry than we need.
We could deform an $N=(4,4)$ theory to e.g. an $N=(0,4)$ heterotic (see e.g. \cite{Wi2} for
an introduction to similar theories). These theories are geometrical, corresponding to
bundles over K3, and are also related to chiral de Rham. They possess the desired elliptic genus $2\phi_{0,1}$. The moduli space of these theories is 90-dimensional, far larger
than that of $N=(4,4)$, so there is a much greater chance for some larger symmetry groups.
However, the breaking of the left-moving $N=4$ supersymmetry seems to destroy the justification
for decomposing the elliptic genus into $N=4$ superconformal characters --- a step
 crucial to Mathieu moonshine. It is tempting to partially break the supersymmetry, say to
 $N=(4,1)$, but such theories seem to possess the full $N=(4,4)$ supersymmetry. (I thank
 Ilarion Melnikov for discussions involving this point.)  

The difficulty in interpreting Mathieu Moonshine in terms of
K3 sigma models (e.g. the interesting work of Taormina-Wendland \cite{TW} is still far from
realising $M_{24}$) leads one to consider the unhappy possibility that the connection with K3 is perhaps accidental.
 The relevant (twisted) elliptic genera are so heavily constrained that
there are bound to be empty coincidences. For example, any of the 71 or so $c=24$
holomorphic rational VOA --- e.g. one associated to the Leech lattice --- will have a
VOA character (a.k.a. partition function) equal to $J(\tau)+c$ for some constant $c\in\bbZ_{\ge 0}$,
again because it is so severely constrained. The coefficients of $J$
(with or without $c$) will have an interpretation as dimensions of Monster representations as
we know,
but conjecturally only one of those VOAs actually carries a nontrivial action of the Monster $\mathbb{M}$. It takes a
(slightly) deeper analysis to rule out Monster actions on these other VOAs. Could this
relation of K3 sigma models to our Mathieu Moonshine be likewise illusory? After all, it is
clear that the Jacobi forms of Umbral Moonshine cannot have a direct interpretation as elliptic genera,
when the group is not $M_{24}$.

Similarly, perhaps we shouldn't regard $M_{24}$ as sacrosanct. The evenness property of Theorem B hints perhaps that the symmetry is somewhat larger. The
analysis of \cite{GHV3} and our Theorem C hints that the `ultimate' symmetry lies
somewhere between $M_{24}$ and Conway. The split extension
$\bbZ_2^{12}\sdprod M_{24}$, a maximal subgroup of Co$_0$, is the only such group
and seems worth a look. (The evenness property was used in \cite{CHM} to prove
one of the Umbral Moonshine conjectures, and together with Theorem A proves
that the elliptic genus of Enriques surfaces  decomposes into a sum of $M_{12}$ characters.)

The appearance of $\Gamma_0(|g|)$ in Lemma 1 is what one would expect from  the CFT orbifold story.
 Let us quickly review that basic theory, as developed by \cite{DVVV},\cite{DW},\cite{Ba},
and others. Let $\mathcal{V}$ be a (bosonic) rational
VOA with automorphisms containing some finite group $G$.
By $\cV^G$ we mean the vertex operator subalgebra
consisting of all fixed-points of $G$ in $\cV$. Conjecturally, $\cV^G$  is also rational; in the
mathematical literature it is called the \textit{orbifold} of $\cV$ by $G$ (the orbifold construction
means something a little different in the physics literature).

We are interested here in the simplest case, where $\cV$ only has a single irreducible module
(namely itself). The twisted modules $M(g)$ of
 $\cV$ are  parametrised by (a subset of the) conjugacy class representatives
$g$, and (a subgroup of) the centraliser of $g$ in $G$  acts on them, so we can define \textit{twisted twining
characters} $Z_{g,h}(\tau)=\mathrm{Tr}_{M_g}(h q^{L_0-c/24})$. The special case
$Z_{1,h}(\tau)$ are sometimes called \textit{McKay-Thompson series} in the mathematics
literature. Then for any $\left(\begin{matrix} {a\atop c}{b\atop d}\end{matrix}\right)\in\mathrm{SL}_2(\bbZ)$, $Z_{g,h}\left(\frac{a\tau+b}{c\tau+d}\right)$ will equal
$Z_{g^ah^c,g^bh^d}(\tau)$ up to some phase. For  general $g,h$, $Z_{g,h}(\tau)$ will thus be
a modular function with multiplier for the group of all matrices  $\left(\begin{matrix} {a\atop c}{b\atop d}\end{matrix}\right)\in\Gamma(\mathrm{lcm}(|g|,|h|))$. However, $\left(\begin{matrix} {a\atop c}{b\atop d}\end{matrix}\right)\in\Gamma_0(|h|)$
will send the McKay-Thompson series $Z_{1,h}$ to another, $Z_{1,h^d}$ (up to a phase).
The integer $d$ will be coprime to $|h|$, and sends a character value of $h$ to a Galois
associate. When the McKay-Thompson series have integer coefficients (e.g. in Monstrous
or Mathieu Moonshine), these
Galois associates have the same value, and  $Z_{1,h^d}=Z_{1,h}$. In other words, since
we have that integrality, our McKay-Thompson series $Z_{1,h}$ will be modular functions 
(with multiplier) for $\Gamma_0 (|h|)$.

To complete the story, we need to consider equivariant cohomology. Let $B_G$ be a classifying 
space, then the group cohomology $H^*_G(1;N)$ (often denoted $H^*(G;N)$)
with values in a $G$-module $N$ is defined by $H^*(B_G;N)$. Just as $H_G^2(1;U(1))$ 
controls the projective representatives of $G$, and $\alpha\in H^3_G(G;\bbZ)$ ($G$ acting on
itself by conjugation) parametrises the possible orbifold fusion rings (all given by
twisted equivariant $K$-theory $^\alpha K^0_G(G)$), the representation theory of $\cV^G$
is controlled by a 3-cocycle $\alpha\in H^3_G(1;U(1))$. In particular, when $\alpha$ is trivial,
the modules $M$ of $\cV^G$ conjecturally are in natural one-to-one correspondence with pairs $(g,\psi)$,
where $g$ is a conjugacy class representative in $G$ and $\psi$ is an irrep of the centraliser 
$C_G(g)$. Their  characters $\chi_{M}(\tau)=\mathrm{Tr}_Mq^{L_0-c/24}$ form a vector-valued modular function for
SL$_2(\bbZ)$ with multiplier explicitly defined in terms of $G$ \cite{DVVV}. The twisted modules $M(g)$ of
 $\cV$ are also (conjecturally) parametrised by all conjugacy class representatives
$g$, and the full centraliser of $g$ in $G$  acts on them Then $Z_{g^k,h^k}(\tau)=Z_{g,h}(\tau)$ and $Z_{g,h}\left(\frac{a\tau+b}{c\tau+d}\right)
=Z_{g^ah^c,g^bh^d}(\tau)$. When the cocycle is not trivial --- the generic case --- a subset of
these pairs $(g,\psi)$ and $(g,h)$ will not be viable. Again, the characters of $\cV^G$
should yield a vector-valued modular form, but the multiplier is more complicated (though
known \cite{CGR}). Similarly, the relations $Z_{g^k,h^k}(\tau)=Z_{g,h}(\tau)$ and $Z_{g,h}\left(\frac{a\tau+b}{c\tau+d}\right)
=Z_{g^ah^c,g^bh^d}(\tau)$ have to be adjusted by phases determined by $\alpha$.

We have order-12 phases in Mathieu Moonshine (e.g. the $e^{2\pi \I cd/(|g|h)}$ in \eqref{phin}), so we would expect (supposing the
orbifold theory extends to our $N=4$ setting) the relevant 3-cocycle $\alpha$ to have order 12.
Indeed, \cite{SE} computed that $H^3_{M24}(1;U(1))\cong \bbZ_{12}$, and \cite{Ron}
explicitly verified that a generating cocycle yields the phases appearing in Mathieu Moonshine. Moreover, 
$H^3_{M23}(1;U(1))\cong 1$ \cite{Mil}, and so the phases for elements lying in $M_{23}$ should
be trivial, and this indeed is what is observed. It is also observed that many pairs $(g,h)$
are not viable, as expected since $\alpha$ is nontrivial. If some group $G$ larger than $M_{24}$
(e.g. Co$_0$ or $2^{12}.M_{24}$) actually acts here, then for this reason we'd require $H^3_G(1;U(1))$ to
contain an order-24 element. (The cohomology of Co$_0$ has not been computed yet).
It is observed in \cite{Ron} that some twisted twining elliptic genera vanish even when there is
no cohomological obstruction --- this would be surprising for a bosonic theory but elliptic
genera (being a signed trace) often vanishes in nontrivial theories so this is not mysterious here.

In summary, the perfect formal fit of twisted twining elliptic genera here with the general
theory of VOA orbifolds, lends strong support to the belief that there is a vertex operator
superalgebra underlying Mathieu Moonshine.

In the case of Monstrous Moonshine, the phases are 24th roots of 1, so we would expect 
the orbifold
$(\cV^\natural)^{\mathbb{M}}$ to be governed by an \textit{order 24} cocycle. We would also expect
this $\alpha$ to obstruct certain pairs $Z_{g,h}$. Indeed, empirical observations encapsulated
in Norton's \textit{Generalised Monstrous Moonshine} \cite{No} make exactly this point:
whether a pair is viable or not is captured by Norton's Fricke dichotomy. This should be equivalent
to the cohomology condition. This relation between cohomology and Norton's mysterious
Fricke condition seems to be new. It would be very interesting to make that relation more
explicit.

                    \newcommand\bibx[4]   {\bibitem{#1} {#2:} {\sl #3} {\rm #4}}

\vspace{0.2cm}\addtolength{\baselineskip}{-2pt}
\begin{footnotesize}
\noindent{\it Acknowledgements.} It is my pleasure to thank
Matthias Gaberdiel and, more generally, the conference participants of the
Mathieu Moonshine workshop at ETH in July 2011, for discussions.
I thank Thomas Creutzig, Matthias Gaberdiel, Gerald H\"ohn, Ilarion Melnikov, Roberto
Volpato and Katrin Wendland for comments involving an earlier draft of this manuscript.

\end{footnotesize}

\end{document}